\documentclass[a4paper, 10pt, final]{amsart}

\usepackage{xcolor}

\usepackage[bookmarks]{hyperref}
\hypersetup{
    allcolors=black,
    pdftitle={Sharelatex Example},
    bookmarks=true,
    pdfpagemode=FullScreen,
}
\hypersetup{
    pdftitle={Your title here},
    pdfauthor={Your name here},
    pdfsubject={Your subject here},
    pdfkeywords={keyword1, keyword2},
    bookmarksnumbered=true,     
    bookmarksopen=true,         
    bookmarksopenlevel=1,       
    colorlinks=true,            
    pdfstartview=Fit,           
    pdfpagemode=UseOutlines,    
    pdfpagelayout=TwoPageRight
}

\usepackage{amsthm}

\theoremstyle{plain}

\usepackage[latin1]{inputenc}
\usepackage[T1]{fontenc}
\usepackage[english]{babel}
\usepackage{amsmath}
\usepackage{lmodern}

\usepackage[OT2,T1]{fontenc}

\usepackage{amsmath, amsthm, amssymb, amsfonts}

\DeclareMathOperator{\weight}{weight}

\DeclareMathOperator{\I}{I}

\DeclareMathOperator{\Ad}{Ad}

\DeclareMathOperator{\depth}{depth}

\DeclareMathOperator{\har}{har}

\DeclareMathOperator{\Li}{Li}

\theoremstyle{definition}

\newtheorem{Theoreme}{Theoreme}[section]
\newtheorem{Proposition}[Theoreme]{Proposition}
\newtheorem{Theorem}[Theoreme]{Theorem}
\newtheorem{Proposition-Definition}[Theoreme]{Proposition-Definition}
\newtheorem{Lemma-Notation}[Theoreme]{Lemma-Notation}
\newtheorem{Lemma-Definition}[Theoreme]{Lemma-Definition}

\newtheorem{Fact}[Theoreme]{Fact}
\newtheorem{Nota Bene}[Theoreme]{Nota Bene}
\newtheorem{Corollary}[Theoreme]{Corollary}
\theoremstyle{definition}

\newtheorem{Question}[Theoreme]{Question}

\newtheorem{Lemma}[Theoreme]{Lemma}
\newtheorem{Definition}[Theoreme]{Definition}
\newtheorem{Remark}[Theoreme]{Remark}

\newtheorem{Example}[Theoreme]{Example}

\newtheorem{Sub-lemma}[Theoreme]{Sub-lemma}

\input cyracc.def
\DeclareFontFamily{U}{russian}{}
\DeclareFontShape{U}{russian}{m}{n}
        { <5><6> wncyr5
        <7><8><9> wncyr7
        <10><10.95><12><14.4><17.28><20.74><24.88> wncyr10 }{}
\DeclareSymbolFont{Russian}{U}{russian}{m}{n}
\DeclareSymbolFontAlphabet{\mathcyr}{Russian}
\makeatletter
\let\@math@cyr\mathcyr
\renewcommand{\mathcyr}[1]{\@math@cyr{\cyracc #1}}
\makeatother
\newcommand{\sh}{\mathcyr{sh}} 

\makeatletter
\newcounter{subsubsubsection}[subsubsection]
\renewcommand\thesubsubsubsection{\thesubsubsection .\@alph\c@subsubsubsection}
\newcommand\subsubsubsection{\@startsection{subsubsubsection}{4}{\z@}%
                                     {-3.25ex\@plus -1ex \@minus -.2ex}%
                                     {1.5ex \@plus .2ex}%
                                     {\normalfont\normalsize\bfseries}}
\newcommand*\l@subsubsubsection{\@dottedtocline{3}{10.0em}{4.1em}}
\newcommand*{\subsubsubsectionmark}[1]{}
\makeatother

\usepackage[top=2.5cm, bottom=2.5cm, left=2.5cm, right=2.5cm]{geometry}

\author{David Jarossay}

\title{Non-vanishing of certain cyclotomic multiple harmonic sums and application to the non-vanishing of certain $p$-adic cyclotomic multiple zeta values}

\renewcommand{\shorttitle}{Non-vanishing of certain cyclotomic multiple harmonic sums and of certain $p$CMZV's}
 
\address{Ben Gurion university of the Negev, Be'er Sheva`, Israel}

\email{jarossay@post.bgu.ac.il}

\begin{document}

\maketitle

\begin{abstract} We define and apply a method to study the non-vanishing of $p$-adic cyclotomic multiple zeta values. We prove the non-vanishing of certain cyclotomic multiple harmonic sums, and, via a formula proved in another paper, which expresses a cyclotomic multiple harmonic sums as an infinite sum of products of $p$-adic cyclotomic multiple zeta values, this implies the non-vanishing of certain $p$-adic cyclotomic multiple zeta values.
\newline This is part III-3 of \emph{$p$-adic cyclotomic multiple zeta values and $p$-adic pro-unipotent harmonic actions}.
\end{abstract}

\tableofcontents

\numberwithin{equation}{section}

\section{Introduction}

In this paper, $\mathbb{N}$ denotes the set of nonnegative integers and $\mathbb{N}^{\ast}$ the set of positive integers.

\subsection{$p$-adic zeta values and non-vanishing}

The values of the Riemann zeta function at positive integers, 
$\zeta(n) = \sum\limits_{m\geq 1} \frac{1}{m^{n}}$ for $n \in \mathbb{N}_{\geq 2}$, have $p$-adic analogues : the values $L_{p}(n,\omega^{1-n})$ where $L_{p}$ is the $p$-adic the Kubota-Leopoldt zeta function, and $\omega$ is Teichmüller's character. To our knowledge, there are only few methods to study the non-vanishing of these values, and few results on it. 

A non-trivial theorem of Soulé \cite{Soule} says that, $L_{p}(n,\omega^{1-n})\not=0$ if $n$ is odd and $p$ is a regular prime number. It is easy to prove that $L_{p}(n,\omega^{1-n})=0$ when $n$ is even, and that $L_{p}(n,\omega^{1-n})\not=0$ if $n$ is odd such that $p-1 | n-1$. We also have the following equation  \cite{Washington} :
\begin{equation} \label{eq:inversion N=1 d=1} 
\sum_{0<m<p} \frac{1}{m^{n}} = \sum_{l\geq 0} {-n \choose l} p^{l} L_{p}(n+l,\omega^{1-n-l}) .
\end{equation}
This has the following known consequence : given that $\displaystyle\sum_{0<m<p} \frac{1}{m^{n}}$ is a strictly positive real number, and thus is non-zero, by (\ref{eq:inversion N=1 d=1}), there exists $l\in \mathbb{N}^{\ast}$ such that $L_{p}(n+l,\omega^{1-n-l}) \not= 0$. Moreover, since this is true for all $n$, there are infinitely many positive integers $n'$ such that $L_{p}(n',\omega^{1-n'}) \not= 0$. 

In view of what follows, let us formulate it differently : for all $n_{0} \geq 1$, we have

\begin{equation} \label{eq:non-vanishing Lp}
(L_{p}(n,\omega^{1-n}))_{n\geq n_{0}} \not= 0 ,
\end{equation}

i.e. the sequence $(L_{p}(n,\omega^{1-n}))_{n\geq n_{0}}$ is not the zero sequence in $\mathbb{Q}_{p}^{\mathbb{N}}$.

\subsection{A generalization : $p$-adic cyclotomic multiple zeta values}

A generalization of the zeta values are the cyclotomic multiple zeta values (CMZV's). For all $d \in \mathbb{N}^{\ast}$ and $(n_{1},\ldots,n_{d}) \in (\mathbb{N}^{\ast})^{d}$ and $\xi_{1},\ldots,\xi_{d}$ roots of unity such that $(n_{d},\xi_{d}) \not= (1,1)$ :

\begin{equation} \label{eq:serie CMZV} \zeta((n_{i})_{d};(\xi_{i})_{d}) = \sum_{0<m_{1}<\ldots<m_{d}} \frac{ (\frac{\xi_{2}}{\xi_{1}})^{m_{1}} \ldots (\frac{1}{\xi_{d}})^{m_{d}}}{m_{1}^{n_{1}} \ldots m_{d}^{n_{d}}} .
\end{equation}

They have an expression as iterated integrals : for $(x_{n},\ldots,x_{1})=(\underbrace{0,\ldots,0}_{n_{d}-1},\xi_{d},\ldots,\underbrace{0,\ldots,0}_{n_{1}-1},\xi_{1})$,

\begin{equation} \label{eq:iterated integral CMZV}
\zeta((n_{i})_{d};(\xi_{i})_{d}) = \int_{0}^{1} \frac{dt_{n}}{t_{n}-x_{n}}  \int_{0}^{t_{n}} \ldots  \int_{0}^{t_{3}} \frac{dt_{2}}{t_{2}-x_{2}} \int_{0}^{t_{2}} \frac{dt_{1}}{t_{1}-x_{1}} .
\end{equation}

$p$-adic cyclotomic multiple zeta values ($p$CMZV's) are defined as $p$-adic analogues of the iterated integrals (\ref{eq:iterated integral CMZV}), where $p\nmid N$. 

In our most general definition \cite{I-1, I-3}, which generalizes the definitions in \cite{Deligne Goncharov, Furusho 1, Furusho 2, Unver 1, Unver 2, Yamashita}, they depend on a parameter $\alpha \in \mathbb{Z} \cup \{\pm \infty\} - \{0\}$ (it is the number of iterations of the crystalline Frobenius of the pro-unipotent fundamental groupoid of $\mathbb{P}^{1} - \{0,\mu_{N},\infty\}$, in the sense of \cite{Deligne}.)

They are elements $\zeta_{p,\alpha} \big((n_{i})_{d};(\xi_{i})_{d}\big)$ of the extension $K_{p}$ of $\mathbb{Q}_{p}$ generated by the $N$-th roots of unity in $\overline{\mathbb{Q}_{p}}$, where $d$, the $n_{i}$'s and the $\xi_{i}$'s are as above (we do not require $(n_{d},\xi_{d}) \not= (1,1)$). And, by \cite{Coleman}, equation (4) p. 173, and by \cite{I-3} corollary 2.2.2, we have
 
\begin{equation} \label{eq:zeta and L} \zeta_{p,\alpha}(n;1) \in \mathbb{Q}^{\times} L_{p}(n,\omega^{1-n}) .
\end{equation}

\subsection{The computation of $p$-adic cyclotomic multiple zeta values and its consequences}

In \cite{I-1, I-2, I-3}, we have computed the $p$CMZV's as sums of series, i.e. we have found a $p$-adic analogue of equation (\ref{eq:serie CMZV}).

The cyclotomic multiple harmonic sums are the following explicit algebraic numbers, where $m$, $d$ and the $n_{i}$'s $(1\leqslant i \leqslant d)$ are positive integers, and the $\xi_{i}$'s $(1\leqslant i \leqslant d+1)$ are $N$-th roots of unity :

\begin{equation} \label{eq:multiple harmonic sum bis} \tilde{h}_{m_{0},m}\big((n_{i})_{d};(\xi_{i})_{d+1}\big) = \displaystyle \sum_{m_{0}<m_{1}<\ldots<m_{d}<m} \frac{\xi_{1}^{m_{0}} \big( \frac{\xi_{2}}{\xi_{1}} \big)^{m_{1}} \cdots \big( \frac{\xi_{d+1}}{\xi_{d}}\big)^{m_{d}}\big(\frac{1}{\xi_{d+1}}\big)^{m}}{m_{1}^{n_{1}} \ldots m_{d}^{n_{d}}} ,
\end{equation}

and we denote by $\tilde{h}_{m} = \tilde{h}_{0,m}$.

Our results of \cite{I-2} include the following equation, which is a large generalisation of equation (\ref{eq:inversion N=1 d=1}) via equation (\ref{eq:zeta and L}) (the left-hand side below is an absolutely congvergent series in $K_{p}$) :

\begin{multline} \label{eq:theequation}
\sum_{d'=0}^{d} \sum_{l_{d'},\ldots,l_{d}\geq 0} \bigg\{ \prod_{i=d'+1}^{d}{-l_{i} \choose n_{i}} \bigg\} \xi_{d'}^{-p^{\alpha}}
\zeta_{p,\alpha}\big((n_{d+1-i}+l_{d+1-i})_{d-d'};(\xi_{d-i})_{d-d'} \big) \zeta_{p,\alpha}\big((n_{i})_{d'};(\xi_{i})_{d'} \big)
\\ = (p^{\alpha})^{\sum_{i=1}^{d} n_{i}} \tilde{h}_{p^{\alpha}}\big((n_{i})_{d};(\xi_{i})_{d+1}\big) .
\end{multline}

In \cite{II-1} and other papers, we use our computation to relate the algebraic theory of $p$CMZV's and the explicit formulas. To do this, we have the following paradigm :

1) We introduce the adjoint $p$CMZV's \cite{II-1}, defined as follows 
\begin{multline} \label{eq:adjoint multiple zeta values} \zeta_{p,\alpha}^{\Ad}\big((n_{i})_{d};(\xi_{i})_{d};l\big) = 
\\ \Big( \sum_{d'=0}^{d} \sum_{\substack{l_{d'},\ldots,l_{d}\geq 0 \\ l_{d'}+\ldots+l_{d}= l}} \bigg\{ \prod_{i=d'+1}^{d}{-l_{i} \choose n_{i}} \bigg\} \xi_{d'}^{-p^{\alpha}}
\zeta_{p,\alpha}\big((n_{d+1-i}+l_{d+1-i})_{d-d'};(\xi_{d-i})_{d-d'} \big) \zeta_{p,\alpha}\big((n_{i})_{d'};(\xi_{i})_{d'} \big)\Big)_{l \in \mathbb{N}} ,
\end{multline}

(the term ``adjoint'' refers to an adjoint action on the pro-unipotent fundamental groupoid of $\mathbb{P}^{1} - \{0,\mu_{N},\infty\}$). For $\Lambda$ a formal variable, we also introduce the $\Lambda$-adjoint $p$CMZVs as 

\begin{equation} \label{eq:adjoint multiple zeta values} \zeta_{p,\alpha}^{\Lambda,\Ad}\big((n_{i})_{d};(\xi_{i})_{d+1}\big) = \sum_{l=0}^{\infty} \Lambda^{l}
\zeta_{p,\alpha}^{\Ad}\big((n_{i})_{d};(\xi_{i})_{d+1};l\big) . 
\end{equation}

Indeed, what appears in all our formulas in \cite{I-1, I-2, I-3} are the adjoint $p$CMZV's, and not the $p$CMZV's themselves (however, $p$CMZV's and adjoint $p$CMZV's can be expressed easily in terms of each other : \cite{II-1} corollary 1.2.5.)
In our terminology, via equation via equation (\ref{eq:zeta and L}), equation (\ref{eq:non-vanishing Lp}) is equivalent to : for any $\alpha$ and $n$, $\zeta_{p,\alpha}^{\Lambda\Ad}(n)\not=0$; and
equation (\ref{eq:theequation}) is reformulated as 
\begin{equation} \label{eq:theequation adjoint} \sum_{l=0}^{\infty}
\zeta_{p,\alpha}^{\Ad}\big((n_{i})_{d};(\xi_{i})_{d+1};l\big) = (p^{\alpha})^{\sum_{i=1}^{d} n_{i}} \tilde{h}_{p^{\alpha}}\big((n_{i})_{d};(\xi_{i})_{d+1}\big) .
\end{equation}

2) Second part of the paradigm. For any question on $p$CMZV's that we want to tackle using explicit formulas, find first an ``adjoint'' variant of the question for adjoint $p$CMZV's and solve it.

This paradigm has proven to be natural and useful for several questions. Each time, we observe the same thing : what is accessible and natural is not the question for the $p$CMZV's themselves but its ``adjoint'' variant which has to be determined.

\subsection{A method for the non-vanishing of $p$-adic cyclotomic multiple zeta values}

When does a $p$CMZV vanish ? To our knowledge, there is no known general approach to this question. Of course, since we have explicit formulas for $p$CMZV's we can check on each example whether a $p$CMZV vanishes or not. But what we want is to prove the non-vanishing of some infinite families of $p$CMZV's.

In this paper, we apply the paradigm explained above (\S0.3) to this question. This means that we are going to study not the non-vanishing of a $\zeta_{p,\alpha}((n_{i})_{d};(\xi_{i})_{d})$, but the non-vanishing of a 
$\zeta^{\Lambda\Ad}_{p,\alpha}((n_{i})_{d};(\xi_{i})_{d+1})$. Once again, as in several previous papers, we are going to observe that considering the ``$\Lambda$-adjoint'' variant of a question makes it accessible.

Let us already see the $N=1$ case. In this case, equation (\ref{eq:theequation}) is simplified into

\begin{equation} \label{eq:theequationN=1} 
\sum_{d'=0}^{d} \sum_{l_{d'},\ldots,l_{d}\geq 0} \bigg\{ \prod_{i=d'+1}^{d}{-l_{i} \choose n_{i}} \bigg\}
\zeta_{p,\alpha}\big((n_{d+1-i}+l_{d+1-i})_{d-d'}\big) \zeta_{p,\alpha}\big((n_{i})_{d'}\big)
\\ = (p^{\alpha})^{\sum_{i=1}^{d} n_{i}} h_{p^{\alpha}}\big((n_{i})_{d}\big) ,
\end{equation}
i.e. 
\begin{equation} \label{eq:theequationN=1 adjoint}
\sum_{l=0}^{\infty}\zeta_{p,\alpha}^{\Ad}\big((n_{i})_{d};l\big) =
 (p^{\alpha})^{\sum_{i=1}^{d} n_{i}} h_{p^{\alpha}}\big((n_{i})_{d}\big) .
\end{equation}

For $p^{\alpha}>d$, by (\ref{eq:multiple harmonic sum bis}), the right-hand side in (\ref{eq:theequationN=1}) is a strictly positive real number, thus is non-zero. By (\ref{eq:theequationN=1 adjoint}), we deduce

\begin{equation} \label{eq:non-vanishing N=1} \zeta_{p,\alpha}^{\Lambda \Ad}((n_{i})_{d}) \not=0 ,
\end{equation}

i.e. there exists $l \in \mathbb{N}_{\geq 0}$ such that $\zeta_{p,\alpha}^{\Ad} ((n_{i})_{d};l) \not=0$. We regard this fact (\ref{eq:non-vanishing N=1}) as the $p$-adic analogue of the following basic fact : all multiple zeta values
$\displaystyle\zeta((n_{i})_{d}) = \sum_{0<m_{1}<\cdots<m_{d}} \frac{1}{m_{1}^{n_{1}} \cdots m_{d}^{n_{d}}}$ are strictly positive real numbers and thus are non-zero.

More generally, for any $N$, by equation (\ref{eq:theequation adjoint}), $\tilde{h}_{p^{\alpha}}((n_{i})_{d};(\xi_{i})_{d+1}) \not= 0$ implies the following equivalent statements: 

(a) $\zeta_{p,\alpha}^{\Lambda \Ad} ((n_{i})_{d};(\xi_{i})_{d+1}) \not=0$

(b) there exists $l \in \mathbb{N}_{\geq 0}$ such that $\zeta^{\Lambda \Ad} ((n_{i})_{d};(\xi_{i})_{d+1};l) \not=0$.

(c) There exists $0 \leq d'\leq d$ and $l_{d'},\ldots,l_{d} \geq 0$ such that  
$\zeta_{p,\alpha}\big((n_{d+1-i}+l_{d+1-i})_{d-d'};(\xi_{d-i})_{d-d'}\big)\not=0$ or $\zeta_{p,\alpha}\big((n_{i})_{d'};(\xi_{i})_{d'} \big)\not=0$ 

Thus, instead of studying the non-vanishing of a $\zeta_{p,\alpha}((n_{i})_{d};(\xi_{i})_{d})$, we are going to the non-vanishing of $\tilde{h}_{p^{\alpha}}((n_{i})_{d};(\xi_{i})_{d+1})$; actually, we are going to study more generally the non-vanishing of cyclotomic multiple harmonic sums (\ref{eq:multiple harmonic sum}). By the above discussion, this is a way to study the ``$\Lambda$-adjoint'' variant of the question of the non-vanishing. (In the end of the paper (\S5) we will discuss the relation between the non-vanishing of $\Lambda$-adjoint $p$CMZV's and the non-vanishing of adjoint $p$CMZV's.)

An advantage of this approach is that $h_{p^{\alpha}}((n_{i})_{d};(\xi_{i})_{d})$ is an algebraic number and is given by a simple formula. We are going to see that in certain cases we can prove that it indeed does not vanish. By contrast, any (adjoint) $p$-adic cyclotomic multiple zeta value is conjecturally either zero or a transcendental number, and is given by a more complicated formula.

This illustrates our general idea that, in order to tackle questions on $p$CMZV's via explicit formulas, we must consider their $\Lambda$-adjoint variants in order to make them accessible.

This also checks, although in an unexpected way, a prediction of Furusho (\cite{Furusho 1}, Remark 2.20) : ``\emph{The author guesses more generally that problems on $p$-adic MZV's related to $p$-adic transcendental number theory (such as the problem of proving the $p$-adic version $\zeta_{p}(3) \not\in \mathbb{Q}$ of Apéry's result) could be translated into problems in algebraic number theory}''.

\subsection{Results} 

From now on, for the simplicity of the notations, we are going to use a different notation for multiple harmonic sums :

\begin{equation} \label{eq:multiple harmonic sum} h_{m_{0},m}\big((n_{i})_{d};(\xi_{i})_{d}\big) = \displaystyle \sum_{m_{0}<m_{1}<\ldots<m_{d}<m} \frac{\xi_{1}^{m_{1}} \cdots \xi_{d}^{m_{d}}}{m_{1}^{n_{1}} \ldots m_{d}^{n_{d}}} .
\end{equation}

Using various elementary techniques, we are going to find necessary conditions on the parameters, namely, $d$, $m$, $(n_{i})_{d}$, $(\xi_{i})_{d}$, to have $h_{m}\big((n_{i})_{d};(\xi_{i})_{d}\big)=0$ i.e. as soon as these conditions are not satisfied, the multiple harmonic sum is non-zero. The techniques will be : observations on the $p$-adic absolute value of multiple harmonic sums, combined with results on the distribution of prime numbers ; reasonings of elementary field theory applied to cyclotomic fields ; observations on the complex absolute value of multiple harmonic sums ; reasonings on alternating series ; using properties of periods.

The main result is the following.

\begin{Theorem}
(i) (\emph{Conditions on $d$ and $m$}) Let $x_{0}$ and $\Delta(x_{0})$ be real numbers as in \cite{RS, KL} with $x_{0}\geq d+1$. If that $m>x_{0}d$ and $\Delta(x_{0})>d+1$, then all multiple harmonic sums $h_{m}\big((n_{i})_{d};(\xi_{i})_{d}\big)$ are non-zero.

(ii) (\emph{Conditions on $\xi_{i}$'s}) If for all $i$ we have $[\mathbb{Q}(\xi_{i}):\mathbb{Q}] = l_{i}\geq m- m_{0}$ with $l_{i}$ and $l_{i'}$ coprime for all $i\not=i'$, then all multiple harmonic sums $h_{m_{0},m}\big((n_{i})_{d};(\xi_{i})_{d}\big)$ and their subsums (see definition \ref{def 2.2}) are non-zero.

If $[\mathbb{Q}(\xi_{1}):\mathbb{Q}]\geq m- m_{0}$ or if $[\mathbb{Q}(\xi_{2}):\mathbb{Q}]\geq m- m_{0}$, then all harmonic sums $h_{m_{0},m}(n_{1},n_{2};\xi_{1},\xi_{2})$ are non-zero.

(iii) (\emph{Conditions on $n_{i}$'s}) If $n_{d}> \frac{\log\big({m-m_{0}-1 \choose d} - 1\big)}{ \log(\frac{m_{0}+d+1}{m_{0}+d})}$ then all multiple harmonic sums $h_{m_{0},m}\big((n_{i})_{d};(\xi_{i})_{d}\big)$ are non-zero.

If $h_{m_{0}+i+2,m}((n_{j});(\xi_{j})_{i+1\leq j \leq d})$ and
$h_{m_{0}+i+1,m}((n_{j});(\xi_{j})_{i+1\leq j \leq d})$ are non-zero, then, for 
\newline $n_{i}> \frac{\log  \Bigg( \displaystyle\frac{\bigg|h_{m_{0}+i+2,m}((n_{j});(\xi_{j})_{i+1\leq j \leq d})\bigg|}{\bigg| h_{m_{0}+i+1,m}((n_{j});(\xi_{j})_{i+1\leq j \leq d}) \bigg|} \Bigg)}{\log\big( \frac{i+1}{i}\big)}$, $h_{m_{0},m}((n_{i})_{d},(\xi_{i})_{d})$ is non-zero.
\end{Theorem}

By \cite{RS} we have $\Delta(x_{0}) \underset{x_{0} \rightarrow \infty}{\rightarrow} \infty$, as we will explain in \S1. Thus (i) above applies to all values of $d$: it implies that for any $d$, we can find $m_{0}(d)$ such that $m\geq m_{0}(d)\Rightarrow h_{m}\big((n_{i})_{d};(\xi_{i})_{d}\big)\not=0$: this leaves a finite number of values of $m$ to study. Combined with (iii), we also have, for each $d$, only finitely many $n_{d}$ to study. By (ii), in depth $\leq 2$ there are also finitely many $(\xi_{i})_{d}$ to study.

Actually in the paper we prove more general results than (i), (ii) (iii), but with less practical conditions, or without application to $p$CMZV's.

By equation (\ref{eq:theequation}) proved in \cite{I-2}, we have the implication theorem 0.1 $\Rightarrow$ corollary 0.2, which is for us the main interest of theorem 0.1:

\begin{Corollary} Let an index $(n_{i})_{d};(\xi'_{i})_{d}$ and a positive integer $m$, such that at least one of the conditions of the theorem is satisfied and that $m=p^{\alpha}$ is a power of a prime number. Let $(\xi_{i})_{d+1}$, such that such that $\xi_{i} = \frac{\xi'_{i+1}}{\xi'_{i}}$ for all $1 \leq i \leq d$.

Then we have $\zeta_{p,\alpha}^{\Lambda\Ad} \big((n_{i})_{d};(\xi_{i})_{d+1}\big)\not= 0$, 

i.e. there exists $l$ such that $\zeta_{p,\alpha}^{\Ad}\big((n_{i})_{d};(\xi_{i})_{d+1};l\big)\not= 0$

i.e. there exists $d' \in \{0,\ldots,d\}$ and $l_{d'+1},\ldots,l_{d} \in \mathbb{N}$ such that
$\zeta_{p,\alpha}\big((n_{i})_{d'};(\xi_{i})_{d'} \big) \not= 0$ or $\zeta_{p,\alpha}\big((n_{d+1-i}+l_{d+1-i})_{d-d'};(\xi_{d-i})_{d-d'} \big)\not= 0$.
\end{Corollary}

As a conclusion, it is possible to prove at least in certain cases the non-vanishing of a $\Lambda$-adjoint $p$CMZV, and there exists at least one natural strategy to study the non-vanishing of $p$CMZV's.

Remarks :

- as the reader will see, although the theorem 0.1 is formulated for \emph{cyclotomic} multiple harmonic sums, it remains true if the $\xi_{i}$'s are replaced with more general numbers (depending on the part of the theorem we are considering), not necessarily roots of unity. However, for simplicity we omit the details about this, because what matters to us is the application to $p$CMZV's.

- we also deduce results of non-vanishing for cyclotomic multiple harmonic values (\S1.3) and cyclotomic multiple zeta values (\S3.2).

- the hypothesis $p\nmid N$ can be removed by using the results of \cite{III-1}, in which we extend the notion of Ad$p$CMZV's and equation (\ref{eq:theequation adjoint}) to the case where $p|N$. We omit the details about this in the paper for the simplicity of the exposition.

Plan of the paper : the results of the theorem are proved in \S1,\S2,\S3, together with auxiliary results and complements. In \S4 we focus on the $N=2$ case and we use an observation on alternating series. In \S5 we focus on adjoint $p$CMZV's by studying the set $\{l \in \mathbb{N}_{\geq 0} \text{ }| \text{ } \zeta_{p,\alpha}^{\Ad} \big((n_{i})_{d};(\xi_{i})_{d};l\big)\not=0\}$, and we also deduce the non-vanishing of cyclotomic multiple harmonic values introduced in \cite{II-1}.

\emph{Acknowledgments.} I thank Daniel Barsky for a discussion during which he taught me the application of Washington's formula (\ref{eq:inversion N=1 d=1}) to a non-vanishing result on values of the Kubota-Leopoldt $p$-adic $L$-function, and he observed that the $N=1$ case of equation (\ref{eq:theequation}) gave a generalization of that non-vanishing result, as explained in the beginning of \S0.4. This is what led me to write this paper. This paper has been supported by NCCR SwissMAP at Universit\'{e} de Gen\`{e}ve, Labex IRMIA at Universit\'{e} de Strasbourg, and ISF grant n°87590031 of Ishai Dan-Cohen at Ben Gurion University of the Negev.

\numberwithin{equation}{subsection}

\section{Conditions on (d,m)}

In this section, we are going to combine an observation on  the $p$-adic absolute value of $h_{m}\big( (n_{i})_{d};(\xi_{i})_{d})$ and results on the distribution of prime numbers.

\subsection{An observation on a $p$-adic absolute value}

The idea is the following. We look at the $p$-adic absolute value of a cyclotomic multiple harmonic sum for a prime number $p$ satisfying certain conditions : in certain cases, the $p$-adic absolute value of a term of the cyclotomic multiple harmonic sum is bigger that the $p$-adic absolute value of the sum of all the other terms. It turns out that the sufficient condition for making this work amounts to find a prime number $p$ in a certain interval. Then, we use results on the distribution of prime numbers, on the existence of prime numbers in certain short effective intervals.

\begin{Lemma} \label{prop 2}
Assume that there exists a positive integer $a$ such that the interval $\big[\big( \frac{m}{d+1} \big)^{1/a}, \big( \frac{m}{d} \big)^{1/a} \big)$ contains a prime number $p>d$.

Then for any positive integers $n_{i}$ ($1\leqslant i \leqslant d$), and for any roots of unity $\xi_{i}$ ($1\leqslant i \leqslant d+1$), we have  $\frak{h}_{m}\big((n_{i})_{d};(\xi_{i})_{d+1}\big)\not=0$.
\end{Lemma}

\begin{proof} For a prime number $p$, and a positive integer $a$, we have the equivalency : $\big( \frac{m}{d+1} \big)^{1/a} \leq p < \big( \frac{m}{d} \big)^{1/a} \Leftrightarrow m \in ]p^{a} d, p^{a}(d+1)]$.

Thus, assume that there exists a prime number $p$ and a positive integer $a$ such that $m \in ]p^{a} d, p^{a}(d+1)]$.

Then, the only integers in $\{1,\ldots,m-1\}$ with $p$-adic valuation equal to $a$ are $p^{a}$, $2p^{a}$, $\ldots$, $dp^{a}$, and all the other integers in $\{1,\ldots,m-1\}$ have $p$-adic valuation in $\{0,\ldots,a-1\}$, since $p^{a}(d+1) \leq p^{a+1}$ by $p>d$. As a consequence, in (\ref{eq:multiple harmonic sum}), the term  $(m_{1},\ldots,m_{d})=(p^{a},2p^{a}\ldots,dp^{a})$ has $p$-adic absolute value ${(p^{a})}^{n_{1}+\ldots+n_{d}}$ and all the other terms have $p$-adic absolute value at most ${(p^{a})}^{n_{1}+\ldots+n_{d} - \min_{1\leqslant i \leqslant d} n_{i}}$. Thus $\big|\har_{m}\big((n_{i})_{d};(\xi_{i})_{d+1}\big)\big|_{p} = {(p^{a})}^{(n_{1}+\ldots+n_{d})} \not= 0$.
\end{proof}

\begin{Remark} The length of the interval	
$[( \frac{m}{d+1} \big)^{1/a}, ( \frac{m}{d} )^{1/a} ]$ is 
$(\frac{m}{d+1})^{1/a} ( (1+\frac{1}{d})^{1/a}-1 )$.
This is a decreasing function of $a$. It tends to $\infty$ when $m \gg d$. It is equivalent to 
$\frac{1}{ad} (\frac{m}{d+1})^{1/a}$ when $d \rightarrow \infty$. 

The length of $[( \frac{m}{d+1})^{1/a}, (\frac{m}{d})^{1/a} ]$ divided by the length of $[ 0,(\frac{m}{d})^{1/a} ]$ is $(\frac{1}{d+1})^{1/a} ( (1+\frac{1}{d})^{1/a}-1 )$. This is again a decreasing function of $a$. It is equivalent to 
$\frac{1}{ad} (\frac{1}{d+1})^{1/a}$ when $d \rightarrow \infty$. In particular, it tends to $0$ when $d \rightarrow \infty$.

Heuristically this tells us that the maximal chance to have a prime number in the interval is for $a=1$, for $m \gg d$ (and perhaps $d$ not too large).
\end{Remark}

\subsection{Use of results on the existence of prime numbers in short effective intervals}

We are now going to use results on the existence of prime numbers in short effective intervals, in order to prove the existence of prime numbers in intervals $\big[\big( \frac{m}{d+1} \big)^{1/a}, \big( \frac{m}{d} \big)^{1/a} \big)$ in certain cases. For simplicity we reduce the discussion to $a=1$.

We are going to use the result of \cite{RS} following theorem from (see also \cite{KL}, theorem 1.1) : let $x_{0} \geq 4. 10^{18}$ be a fixed constant and let $x > x_{0}$. Then there exists at least one prime $p$ such that $(1 - \Delta^{-1})x < p < x$, where $\Delta$ is a constant depending on $x_{0}$ and is given numerically in table 2 of \cite{RS}. For example, by \cite{RS}, theorem 3, if $x_{0} = 10 726 905 041$, we have $\Delta = 28314000$.

\begin{Lemma} \label{second lemma prime number}Let $x_{0}$ and $\Delta(x_{0})$ as in the theorem 1.1 of \cite{KL}, with $x_{0}\geq d+1$. Then, if $d+1 < \Delta(x_{0})$ and $x_{0}d < m$, then the interval $[\frac{m}{d+1},\frac{m}{d}[$ contains a prime number $p>d$.
\end{Lemma}

\begin{proof} By the theorem mentioned above, to have a prime number in the interval $[\frac{m}{d+1},\frac{m}{d}[$, it suffices to show an inclusion $]x(1 - \frac{1}{\Delta(x_{0})}),x] \subset [\frac{m}{d+1},\frac{m}{d}[$ with $x \in \mathbb{R}$ such that $x \geq \Delta(x_{0})$.
	
We want to find $x$ such that $\left\{ \begin{array}{l} x \geq x_{0} 
\\ \frac{m}{d+1} < x (1 - \frac{1}{\Delta(x_{0})})
\\ x < \frac{m}{d}
\end{array} \right.$. This is equivalent to $\left\{ \begin{array}{l} x \geq x_{0} 
\\ \frac{m}{d+1} (1 - \frac{1}{\Delta(x_{0})})^{-1} < x  < \frac{m}{d}
\end{array} \right.$.
	
The existence of such a real number $x$ is equivalent to the two equalities $\left\{ \begin{array}{l} \frac{m}{d+1} (1 - \frac{1}{\Delta(x_{0})})^{-1} < \frac{m}{d} 
\\ x_{0} < \frac{m}{d} 
\end{array} \right.$, i.e. 
$\left\{ \begin{array}{l} \frac{d}{d+1} = 1 - \frac{1}{d+1} < 1 - \frac{1}{\Delta(x_{0})}
\\ x_{0}d < m
\end{array} \right.$
i.e. 
$\left\{ \begin{array}{l} d+1 < \Delta(x_{0})
\\ x_{0}d < m
\end{array} \right. .$

In these conditions, a prime number $p$ in the interval $[\frac{m}{d+1},\frac{m}{d}[$ satisfies $p> \frac{m}{d+1} > \frac{x_{0}d}{d+1}$. Thus, for $x_{0}\geq d+1$, we have $p>d$.
\end{proof}

We have the following result \cite{Balliet} : for any positive integers $n$ and $k$ such that $n \geq k$ and $2 \leq k \leq 519$, the interval $]kn,(k+1)n[$ where  contains at least one prime number. We deduce from this result the following :

\begin{Lemma} \label{first lemma prime number} Assume $d<260$ and $m>\max(2d(d+1),520^{2}d,\frac{520d(d+1)}{260 - d})$. then the interval $]\frac{m+1}{d},\frac{m}{d}]$ contains a prime number $p>d$.
\end{Lemma}

\begin{proof} By the theorem cited above, 
it is sufficient to find $k$ and $n$ such that $n \geq k$ and $2 \leq k \leq 519$, the interval $]kn,(k+1)n[$ is included in $] \frac{m}{d+1},\frac{m}{d}]$, and to have $\frac{m}{d+1}>d$, i.e. $m>d(d+1)$.
	
Thus, given $m$ and $d$, we want to find $k$ and $n$ such that 
$\left\{\begin{array}{l} \displaystyle \frac{m}{d+1}\leq kn \\
(k+1)n \leq \frac{m}{d} \displaystyle
\\ n \geq k
\\ 2 \leq k \leq 519
\end{array} \right. $ i.e. $\left\{\begin{array}{l} \displaystyle \frac{m}{n(d+1)}\leq k < \frac{m}{nd} - 1 
\\ n \geq k
\\ 2 \leq k \leq 519
\end{array} \right. $.
	
If we have $|\frac{m}{nd} - 1 - \frac{m}{n(d+1)}| > 1$, i.e. 
$|\frac{m}{nd(d+1)} - 1| > 1$, 
then the interval $]\frac{m}{n(d+1)}, \frac{m}{nd} - 1]$ contains at least one integer, and we can define $k$ to be such an integer. 
	
Thus, let us analyse the condition : $|\frac{m}{nd(d+1)} - 1| > 1$
	
First case : $\frac{m}{nd(d+1)} - 1<0$, i.e. $\frac{m}{nd(d+1)}<1$ i.e. $\frac{m}{d(d+1)} < n$.
	
The condition $|\frac{m}{nd(d+1)}-1| > 1$ is equivalent to 
$1 - \frac{m}{nd(d+1)} > 1$ i.e. $\frac{m}{nd(d+1)} < 0$. This is impossible.
	
Second case : $\frac{m}{nd(d+1)} - 1>0$ i.e. $\frac{m}{nd(d+1)}>1$. 
	
Then the condition is equivalent to $\frac{m}{nd(d+1)}>2$, i.e. $n < \frac{m}{2d(d+1)}$.
	
Let $n$ satisfying this condition. Let $k$ be an integer in $]\frac{m}{n(d+1)}, \frac{m}{nd} - 1]$. 
	
We want to see under which condition $n \geq k$ and $2 \leq k \leq 519$. 
	
Since $k < \frac{m}{nd} - 1$, in order to have $n \geq k$, it is sufficient to have $n \geq \frac{m}{nd} - 1$, i.e. $n+1 \geq \frac{m}{nd}$, i.e $n(n+1) \geq \frac{m}{d}$. Since $n(n+1)>n^{2}$, this is implied by $n \geq \sqrt{\frac{m}{d}}$.
	
Since $\frac{m}{(n+1)d} < k < \frac{m}{nd} - 1$, in order to have $2 \leq k \leq 519$, it is sufficient to have $2 \leq \frac{m}{(n+1)d}$ and $\frac{m}{nd} - 1 \leq 519$, i.e.
$n+1 \leq \frac{m}{2d}$ and $\frac{m}{nd} \leq 520$, 
i.e. $\frac{m}{520 d} \leq n \leq \frac{m}{2d} - 1$. 
	
Thus, to find $k$ and $n$ as we want, it suffices to find $n$ such that 
	
$$ \max (\frac{m}{520d}, \sqrt{\frac{m}{d}}) \leq n < \min ( \frac{m}{2d(d+1)},\frac{m}{2d} - 1) $$
	
i.e. the interval $[\max (\frac{m}{520d}, \sqrt{\frac{m}{d}}), \min ( \frac{m}{2d(d+1)},\frac{m}{2d} - 1)[$ must contain at least one integer. 
	
It is sufficient to have $\min ( \frac{m}{2d(d+1)},\frac{m}{2d} - 1) - \max (\frac{m}{520d}, \sqrt{\frac{m}{d}})> 1$.
	
For $\frac{m}{2(d+1)}>1$, we have $\frac{m}{2d(d+1)}\leq \frac{m}{2d} - 1$ and thus $\min ( \frac{m}{2d(d+1)},\frac{m}{2d} - 1)=\frac{m}{2d(d+1)}$
	
For $\sqrt{m} \geq 520 \sqrt{d}$, we have $\frac{m}{520d} \geq \sqrt{\frac{m}{d}}$, and $\max (\frac{m}{520d}, \sqrt{\frac{m}{d}})=\frac{m}{520d}$. 
	
So, for $m>\max(2d(d+1),520^{2}d)$, the condition $\min ( \frac{m}{2d(d+1)},\frac{m}{2d} - 1) - \max (\frac{m}{520d}, \sqrt{\frac{m}{d}})> 1$ becomes 
	
$\frac{m}{2d(d+1)}- \frac{m}{520d} > 1$, i.e. $m(260 - d)
>520d(d+1)$. 

This is equivalent to $260 >d$ and $m>\frac{520d(d+1)}{260 - d}$. 
	
So in the end the conditions are satisfied for $d<260$ and $m>\max(2d(d+1),520^{2}d,\frac{520d(d+1)}{260 - d})$.
	
\end{proof}

Combining the previous lemmas we obtain :

\begin{Proposition} Let $m,d$ be positive integers with $m>d$. Assume that $m>x_{0}d$ and $d+1 < \Delta(x_{0})$ with $x_{0}$ and $\Delta(x_{0})$ as in \cite{KL} and $x_{0}\geq d+1$, or assume that 
$d<260$ and $m>\max(2d(d+1),520^{2}d,\frac{520d(d+1)}{260 - d})$. Then, for all $(n_{i})_{d}$, and $(\xi_{i})_{d}$, $h_{m}((n_{i})_{d};(\xi_{i})_{d})\not=0$.
\end{Proposition}

\begin{proof} This follows lemma \ref{prop 2} combined with
lemma \ref{second lemma prime number} and lemma \ref{first lemma prime number}.
\end{proof}

We have equation (2) in \cite{RS} :
$$ \max\{ \Delta \geq 1 \text{ }|\text{ }\forall x \geq (T_{0}\log T_{0})^{2},\text{ }\exists p \in ]x(1-\Delta^{-1}),x]\} \gg T_{0} \log(T_{0})^{-1/2} . $$
This indicates that we can take $\Delta(x)$ such that $\Delta(x) \underset{x \rightarrow \infty}{\rightarrow} \infty$, and proposition 1.5 applies to all values of $d$.

\subsection{Remarks}

\subsubsection{The other cases}

Behind this method, there is a more general topic : the study of $h_{m}((n_{i})_{d};(\xi_{i})_{d})$ as a function of $m$ regarded as a $p$-adic integer. We have an ``expansion'' of $h_{m}((n_{i})_{d};(\xi_{i})_{d})$ in terms of the $q$-adic expansion of $m$, where $q$ is the cardinality of the residue field of the $p$-adic field $K_{p}$ defined in \S0.2 (see \cite{I-3}, proposition 3.3.1).

Here, we have focused on the first term of the expansion, and we have considered a favorable case in which we can prove that this first term is non-zero. In general, for any prime number $p$, there exists a unique $a$ such that $m \in ]p^{a},p^{a+1}]$. Thus there exists a unique $d'$ with $1 \leq d' \leq p-1$ and $m \in ]p^{a}d',p^{a}(d'+1)]$. The case we have considered in \S1.1 is $d'=d$. If $d'>d$, we can prove, similarly to lemma \ref{prop 2}, the implication 
$$ v_{p}(h_{d'}((n_{i})_{d};(\xi_{i}^{p^{a}})_{d}) < \underset{1\leq i \leq d}{\min} n_{i} \Rightarrow h_{m}((n_{i})_{d};(\xi_{i})_{d})\not=0 , $$
and a similar implication if $d'<d$. We see that it has a stronger hypothesis : an information on the $p$-adic valuation of a cyclotomic multiple harmonic sum, more than its non-vanishing.
The study of the $p$-adic valuation of  $h_{r}((n_{i})_{d};(\xi_{i}^{p^{a}})_{d})$ with $1\leq r \leq p$ is a deeper and more difficult problem. In \cite{II-3}, we have interpreted this problem in terms of the study of the difference between the slopes of Frobenius and the Hodge filtration on the crystalline pro-unipotent fundamental groupoid of $\mathbb{P}^{1} - \{0,\mu_{N},\infty\}$.

\subsubsection{Application to cyclotomic multiple harmonic values}

We have defined and studied in \cite{II-1} the cyclotomic multiple harmonic values (CMHV's) as follows, where $\mathcal{P}_{N}$ is the set of prime numbers such that $p\nmid N$. We use the notation $K_{p}$ from \S0.2. Let $w=((n_{i})_{d};(\xi_{i})_{d+1})$.

(a) For $p \in \mathcal{P}_{N}$, $\har_{p^{\mathbb{N}}}(w) 
= \big( \har_{p^{\alpha}} (w)\big)_{\alpha \in \mathbb{N}} \in K_{p}^{\mathbb{N}}$ is called a $p$-adic CMHV.

(d) For $\alpha \in \mathbb{N}^{\ast}$, let $\har_{\mathcal{P}^{\alpha}}(w) = \big( \har_{p^{\alpha}} (w)\big)_{p  \in  \mathcal{P}_{N}} \in \prod\limits_{p\in\mathcal{P}_{N}} K_{p}$, is called an adelic CMHV.

(c) Let $\har_{\mathcal{P}_{N}^{\mathbb{N}}}(w) =
\big(\har_{p^{\alpha}} (w)\big)_{(p,\alpha)  \in  \mathcal{P}_{N} \times \mathbb{N}} \in \big( \prod\limits_{p\in\mathcal{P}_{N}} K_{p} \big)^{\mathbb{N}}$ is called a ``$p$-adic$\times$adelic'' CMHV.

CMHV's and ($\Lambda$-)adjoint $p$CMZV's are the central objects in our explicit description of the algebraic theory of $p$CMZV's.

Here, we have proved :

\begin{Corollary} All cyclotomic multiple harmonic values are non-zero. More precisely, for a given cyclotomic multiple harmonic value (of any of the three types above), all the terms of rank large enough are non-zero.
\end{Corollary}

\begin{proof} This is a consequence of proposition 1.5 : given $d$ and $(n_{i})_{d};(\xi_{i})_{d})$, for all $m$ large enough, we have $h_{m}((n_{i})_{d};(\xi_{i})_{d})\not=0$.
\end{proof}

\subsubsection{An easier result}

If we ask only for the existence of infinitely many $m$ such that  $\frak{h}_{m}\big((n_{i})_{d};(\xi_{i})_{d+1}\big)\not=0$ (instead of for all $m$ large enough), the result is much simpler to obtain.

\begin{Proposition} There exist infinitely many $m$'s such that $\frak{h}_{m}\big((n_{i})_{d};(\xi_{i})_{d+1}\big)\not=0$.
\end{Proposition}

\begin{proof} The multiple polylogarithms on $\mathbb{P}^{1} - \{0,\mu_{N},\infty\}$ are the following functions : $\Li(\emptyset)=1$ and, for all $d\geq 0$,

$$ \Li(((n_{i})_{d+1};(\xi_{i})_{d+1}\big)) = \sum_{0<m} \tilde{h}_{m}\big((n_{i})_{d};(\xi_{i})_{d+1}\big)\frac{z^{m}}{m^{n_{d+1}}} . $$

By Chen's theorem \cite{Chen} on iterated integrals, multiple polylogarithms are linearly independent over the algebraic functions on $\mathbb{P}^{1} - \{0,\mu_{N},\infty\}$. In particular, each $\Li(((n_{i})_{d+1};(\xi_{i})_{d+1}\big))$ is not a polynomial, thus there are infinitely many $m$ such that $\tilde{h}_{m}\big((n_{i})_{d};(\xi_{i})_{d+1}\big)\not=0$.
\end{proof}

\section{Conditions on $(\xi_{i})_{d}$}

In this section we want to show that if a multiple harmonic sum is zero, this implies certain conditions on the $\xi_{i}$'s, or more precisely on their degrees as algebraic numbers. We are going to use elementary field theory. If $\xi$ is a root of unity of order $N$, we have $[\mathbb{Q}(\xi) : \mathbb{Q}] = \varphi(N)$ where $\varphi$ is Euler's function. We have $\varphi(N) \rightarrow \infty$ when $N \rightarrow \infty$, thus it is equivalent to bound $\varphi(N)$ and to bound $N$.

\subsection{Generalities}

The starting point of this section is the following observation : 

\begin{Fact} \label{fact 2.1}Given a multiple harmonic sum $h_{m_{0},m}((n_{i})_{d};(\xi_{i})_{d})$, and $1 \leq i \leq d$, let the polynomial 
$$ P_{i}(X) = \sum_{m_{0}<m_{1}<\ldots <m_{d}<m} \frac{\xi_{1}^{m_{1}} \cdots \xi_{i-1}^{m_{i-1}} X^{m_{i}-m_{0}} \xi_{i+1}^{m_{i+1}} \cdots \xi_{d}^{m_{d}}}{m_{1}^{n_{1}} \cdots m_{d}^{n_{d}}} . $$ 
We have $h_{m_{0},m}((n_{i})_{d};(\xi_{i})_{d}) = \xi_{i}^{m_{0}}P_{i}(\xi_{i})$, and $P_{i}(X) \in \mathbb{Q}(\xi_{1},\ldots,\widehat{\xi_{i}},\ldots,\xi_{d})[X]$ where $\widehat{\xi_{i}}$ means that $\xi_{i}$ is omitted.  Moreover, if $P_{i}\not=0$ we have $\deg(P_{i}) \leq m - m_{0}-1$. We deduce that if $P_{i}\not=0$ and $h_{m_{0},m}((n_{i})_{d};(\xi_{i})_{d})=0$, we have $P_{i}(\xi_{i})=0$ with $\deg(P_{i}) \leq m - m_{0}-1$, thus 
$[\mathbb{Q}(\xi_{1},\ldots,\xi_{d}):\mathbb{Q}(\xi_{1},\ldots,\widehat{\xi_{i}},\ldots,\xi_{d}) ]\leq m - m_{0} - 1$. In other terms we have the implication :

$$ P_{i}\not=0\text{ and }[\mathbb{Q}(\xi_{1},\ldots,\xi_{d}):\mathbb{Q}(\xi_{1},\ldots,\widehat{\xi_{i}},\ldots,\xi_{d}) ] \geq m - m_{0} \Rightarrow h_{m_{0},m}((n_{i})_{d};(\xi_{i})_{d})\not=0 . $$
\end{Fact}

We now explain how to transform this fact into a statement which does not require an assumption of the type ``$P_{i} \not=0$'' but involve other types of assumptions instead, by an inductive process.

\begin{Definition} \label{def 2.2} A subsum of a multiple harmonic sum 
$h_{m_{0},m}((n_{i})_{d};(\xi_{d}))$
is a multiple harmonic sum of the form $h_{m'_{0},m'}((n_{j})_{d'};(\xi_{j})_{d'}))$ where $m_{0} \leq m'_{0} < m' \leq m$, $d' \leq d$, and where there exists an interval of integers $[a,b] \subset [1,d]$ with $b-(a-1)=d'$ and such that $(n_{j})_{d'} = (n_{i})_{a \leq i \leq b}$ and $(\xi_{j})_{d'} = (\xi_{i})_{a \leq i \leq b}$.
\end{Definition}

\begin{Proposition} \label{prop general field}Let a multiple harmonic sum $h_{m_{0},m}((n_{i})_{d};(\xi_{i})_{d})$. Assume that for all subsums $h_{m'_{0},m'}((n_{j})_{d'};(\xi_{j})_{d'}))$ we have  have 
$[\mathbb{Q}(\xi_{j_{1}},\ldots,\xi_{j_{d'}}):\mathbb{Q}(\xi_{j_{1}},\ldots,\widehat{\xi_{j_{l}}},\ldots,\xi_{j_{d'}})] \geq m' - m'_{0}$, where $\widehat{\xi_{j_{l}}}$ means that $\xi_{j_{l}}$ is omitted.

Then the multiple harmonic sum $\displaystyle\sum_{m_{0}<m_{1}<\ldots <m_{d}<m} \frac{\xi_{1}^{m_{1}} \cdots \xi_{d}^{m_{d}}}{m_{1}^{n_{1}} \cdots m_{d}^{n_{d}}}$ and all its subsums are non-zero.
\end{Proposition}

\begin{proof} We are going to use induction on $(m - m_{0} - d,d)$, where we call $m - m_{0} - d$ the ``length minus depth'', and where we use the lexicographical order on $\mathbb{N}_{\geq 1}^{2}$.
	
By induction we prove the following statement on subsums : 
	
(a) for all $i$, the polynomial $\displaystyle P_{i}(X) = \sum_{m'_{0}<m_{1}<\ldots <m_{d'}<m'} \frac{\xi_{j_{1}}^{m_{1}} \cdots \xi_{i-1}^{m_{i-1}} X^{m_{i}-m_{0}} \xi_{i+1}^{m_{i+1}} \cdots \xi_{j_{d'}}^{m_{d'}}}{m_{1}^{n_{1}} \cdots m_{d}^{n_{d'}}}$
is non-zero.

(b) the multiple harmonic sum $\displaystyle\sum_{m'_{0}<m_{1}<\ldots <m_{d'}<m'} \frac{\xi_{j_{1}}^{m_{1}} \cdots \xi_{j_{d'}}^{m_{d}}}{m_{1}^{n_{1}} \cdots m_{d}^{n_{d}}}$ is non-zero.
	
For $m' - m'_{0}- d' =1$ (and any $d'$), it follows from the definition that each polynomial $P_{i}$ is a monomial multiplied by a non-zero coefficient, and is thus non-zero. And the corresponding multiple harmonic sum has exactly one term ($m_{i}=m_{0}+i$ for all $i$), which is clearly non-zero.

For $d'=1$ (and any $m' - m'_{0}- d'$), each coefficient of $P_{i}$ is a strictly positive real number, thus $P_{i}$ is non-zero. The hypothesis of the proposition implies $[\mathbb{Q}(\xi_{i}):\mathbb{Q}]\geq m' - m'_{0}$ for all $i$. By fact 2.1, this implies that the corresponding harmonic sum is non-zero.

We assume $m' - m'_{0}- d' \geq 2$ and $d' \geq 2$.

Each coefficient of $P_{i}$ is a product of two multiple harmonic sums (when $i=1$ or $i=d$, one of the two multiple harmonic sums is the trivial one equal to $1$). Let us consider the coefficient of a give degree, $X^{m_{i} - m'_{0}}$. For the two multiple harmonic sum factors, the length minus depth are respectively $m'-m_{i} - (d'-i)$ and $m_{i} - m'_{0} - (i-1)$. The sum of these two lengths minus depths is $m'-m_{i} - d'+ i + m_{i} - m'_{0} - i+ 1 = m' - m'_{0} - d' + 1$. Both of them are $\leq m'-m'_{0} - d'$, and at least one of them is $<m' - m'_{0} - d'$ : otherwise their sum would be $\geq 2 (m'-m'_{0} - d')$, whereas we have $m' - m'_{0} - d' + 1 < 2 (m'-m'_{0} - d')$ since $m' - m'_{0} - d' \geq 2$. If one of them is equal to $m'-m'_{0} - d'$, the corresponding depth is strictly smaller than $d'$.

For both factors, we have (length-depth,depth) $<$ ($m'-m'_{0}-d',d'$). By the hypothesis of the proposition, and by the induction hypothesis, we deduce that both multiple harmonic sums are non-zero, thus the coefficient of degree $m_{i} - m_{0}$ of $P_{i}$ is non-zero, thus $P_{i}$ is non-zero.

If we had $\displaystyle\sum_{m'_{0}<m_{1}<\ldots <m_{d'}<m'} \frac{\xi_{j_{1}}^{m_{1}} \cdots \xi_{j_{d'}}^{m_{d}}}{m_{1}^{n_{1}} \cdots m_{d}^{n_{d}}}=0$, we would have $P_{i}(\xi_{i})=0$ with $P_{i} \not=0$, $P_{i} \in \mathbb{Q}(\xi_{1}\cdots \widehat{\xi_{i}} \cdots \xi_{d})[X]$ of degree $\leq m-m_{0} - 1$, whence a contradiction with the hypotheis.

Thus $\displaystyle\sum_ {m'_{0}<m_{1}<\ldots <m_{d'}<m'} \frac{\xi_{j_{1}}^{m_{1}} \cdots \xi_{j_{d'}}^{m_{d}}}{m_{1}^{n_{1}} \cdots m_{d}^{n_{d}}}\not=0$. Whence the result by induction.
\end{proof}

In the next subsections, we are going to give concrete situations in which the proposition applies.

\subsection{The linearly disjoint case}

We assume all fields to be of characteristic $0$. We have :

\begin{Fact} \label{2.1}For all finite extensions $M$ and $L$ of a field $K$, we have $[ML : M] \leq [L : K]$.
\end{Fact}

Let $L_{1},\ldots,L_{d}$ finite extensions of a field $K$. By multiplicativity of the degrees we have 
$[L_{1}\cdots L_{d} : K ] = [L_{1}\cdots L_{d};L_{1}\cdots L_{d-1}] \cdots [L_{1}L_{2} : L_{1}] [L_{1} : K]$ whence, by the previous fact,

$$ [L_{1}\cdots L_{d} : K ] \leq  [ L_{d}; K] \cdots [L_{2} : K] [L_{1} : K] $$

\begin{Definition} When the above inequality is an equality, we say that $L_{1},\ldots,L_{d}$ are linearly disjoint (as extensions of $K$).
\end{Definition}

\begin{Fact} \label{fact linearly disjoint property} If $L_{1},\ldots,L_{d}$ are linearly disjoint, then for any $1\leq i \leq d$,

(a) $L_{1},\ldots,L_{i}$ are linearly disjoint. 

(b) $L_{i+1},\ldots,L_{d}$ are linearly disjoint. 

(c) $[L_{1} \cdots L_{d}; L_{1} \cdots L_{i}] = [L_{i+1},\ldots,L_{d}:K]$

(d) We have $[L_{i}:K]=[ L_{1}\cdots L_{d}; L_{1}\cdots \widehat{L_{i}} \cdots L_{d}]$ where $\widehat{L_{i}}$ means that $L_{i}$ is omitted in the product. In particular,
$[L_{1}\cdots L_{d} : K ] = \prod_{i=1}^{d}[ L_{i}; K]  = \prod_{i=1}^{d}[ L_{1}\cdots L_{d}; L_{1}\cdots \widehat{L_{i}} \cdots L_{d}]$.
\end{Fact}

\begin{proof} By the multiplicativity of degrees we have 
$[L_{1}\cdots L_{d} : K ] = [L_{1} \cdots L_{d}; L_{1} \cdots L_{i}] [L_{1} \cdots L_{i};K]$.

Moreover by fact \ref{2.1} we have 
$[L_{1} \cdots L_{d}; L_{1} \cdots L_{i}] \leq [L_{i+1}\cdots L_{d}:K]$
and 
$[L_{1} \cdots L_{i};K] \leq \prod_{j=1}^{i}[L_{j}:K]$

Whence: $[L_{1}\cdots L_{d} : K ] = [L_{1} \cdots L_{d}; L_{1} \cdots L_{i}] [L_{1} \cdots L_{i};K] \leq [L_{1} \cdots L_{d}; L_{1} \cdots L_{i}]\prod\limits_{j=1}^{i}[L_{j}:K] \leq \prod_{i=1}^{d} [L_{i}:K]$. 

Thus, if the inequality $[L_{1}\cdots L_{d} : K ] \leq \prod_{i=1}^{d} [L_{i}:K]$ is an equality, then in particular all the intermediate inequalities 
$[L_{1} \cdots L_{d}; L_{1} \cdots L_{i}] [L_{1} \cdots L_{i};K] \leq [L_{1} \cdots L_{d}; L_{1} \cdots L_{i}]\prod_{j=1}^{i}[L_{j}:K]$ and 
$ [L_{1} \cdots L_{d}; L_{1} \cdots L_{i}] \leq [L_{i+1},\ldots,L_{d}:K] \leq \prod_{j=i+1}^{d} [L_{j}:K]$, are equalities. The first one amounts to $[L_{1} \cdots L_{i};K] = \prod\limits_{j=1}^{i}[L_{j}:K]$, which proves (a), and the two other ones prove (b) and (c). 

In particular, we have $[L_{1} \cdots L_{d}; L_{1} \cdots L_{d-1}] = [L_{d}:K]$. Since the roles of the $L_{i}$'s are symmetric, more generally we have 
$[L_{1}\cdots L_{d}; L_{1}\cdots \widehat{L_{i}} \cdots L_{d}] = [L_{i}:K]$ for all $i$. Whence (d).
\end{proof}

\begin{Fact} \label{2.7}Assume that $L_{1},\ldots,L_{d-1}$ are linearly disjoint, let $L=\prod_{i=1}^{d-1}L_{i}$, and assume that $L_{d},L$ are linearly disjoint. Then  $L_{1},\ldots,L_{d}$ are linearly disjoint.
\end{Fact}

\begin{proof} By the multiplicativity of degreees, $[L_{1}\ldots L_{d}:K] = [L_{1}\ldots L_{d}:L][L:K] = [LL_{d}:L][L:K]$.

By the linear disjointness of $L_{1},\ldots,L_{d-1}$ we have $[L:K] = \prod_{i=1}^{d-1} [L_{i}:K]$.

By the linear disjointness of $L,L_{d}$, and by fact 2.6 (c), we have $[LL_{d}:L] = [L_{d}:K]$. 

Whence $[L_{1}\ldots L_{d}:K] = \prod_{i=1}^{d} [L_{i}:K]$.
\end{proof}

\begin{Fact} \label{fact linearly disjoint coprime}For $d \geq 2$, if $L_{i} = \mathbb{Q}(\xi_{i})$ for $1 \leq i \leq d$, with $\xi_{i}$ a root of unity of order $l_{i}$ with for all $i\not= i'$, $l_{i}$ and $l_{i'}$ coprime, then

(a) there exists a root of unity $\xi'$ of order $\prod_{i=1}^{d} l_{i}$ such that $\mathbb{Q}(\xi_{1},\ldots,\xi_{d}) = \mathbb{Q}(\xi')$;

(b) $\mathbb{Q}(\xi_{1}),\ldots,\mathbb{Q}(\xi_{d})$ are linearly disjoint.
\end{Fact}

\begin{proof} By induction on $d$.

Assume $d=2$. We denote $\xi_{1}=\xi$, $\xi_{2}= \xi'$, $l_{1}=l$, $l_{2}=l'$. We can assume $\xi=e^{2i\pi / l}$, $\xi'=e^{2i\pi / l'}$, without changing the extensions $\mathbb{Q}(\xi)$ and $\mathbb{Q}(\xi')$. For any $u,v \in \mathbb{Z}$, we have then $\xi^{u}{\xi'}^{v} = e^{2i\pi u / l} e^{2i\pi v/ l'} = e^{2i\pi \frac{ul'+vl}{ll'}}$. Since $l$ and $l'$ are coprime, we can find $u,v$ such that $ul + vl' = 1$, whence $ e^{2i\pi/ll'} \in \mathbb{Q}(\xi,\xi')$, whence $\mathbb{Q}(e^{2i\pi/ll'}) \subset \mathbb{Q}(\xi,\xi')$ and the converse inclusion is clear. In particular $[\mathbb{Q}(\xi,\xi') : \mathbb{Q}] = \varphi(ll')$. By the multiplicativity of the Euler function $\varphi$, since $l$ and $l'$ are coprime, we have $\varphi(ll')=\varphi(l)\varphi(l')$ i.e. $[\mathbb{Q}(\xi,\xi'):\mathbb{Q}] = [\mathbb{Q}(\xi):\mathbb{Q}] [\mathbb{Q}(\xi'):\mathbb{Q}]$.	 

Now, let us assume the result is known for $d-1$ and let us prove the result for $d$. By the induction hypothesis, we have $\mathbb{Q}(\xi_{1},\ldots,\xi_{d-1}) = \mathbb{Q}(\xi')$, with $\xi'$ a root of unity of order $\prod_{i=1}^{d-1} l_{i}$, and $\mathbb{Q}(\xi_{1}),\ldots,\mathbb{Q}(\xi_{d-1})$ are linearly disjoint. $l_{d}$ is coprime with all $l_{i}$, $1 \leq i \leq d-1$, so it is coprime with $\prod_{i=1}^{d-1} l_{i}$. Thus, the orders of $\xi_{d}$ and $\xi'$ are coprime. By the result for $d=2$, this implies that $\mathbb{Q}(\xi_{d}),\mathbb{Q}(\xi')$ are linearly disjoint, and that there exists $\xi''$, root of unity of order $\prod_{i=1}^{d} l_{i}$, such that $\mathbb{Q}(\xi_{d},\xi') = \mathbb{Q}(\xi'')$.

By fact \ref{2.7}, since $\mathbb{Q}(\xi_{d}),\mathbb{Q}(\xi')$ are linearly disjoint and  $\mathbb{Q}(\xi_{1}),\ldots,\mathbb{Q}(\xi_{d-1})$ are linearly disjoint, we deduce that $\mathbb{Q}(\xi_{1}),\ldots,\mathbb{Q}(\xi_{d})$ are linearly disjoint. Thus $\mathbb{Q}(\xi_{1},\ldots,\xi_{d-1},\xi_{d}) = \mathbb{Q}(\xi',\xi_{d}) = \mathbb{Q}(\xi'')$.
\end{proof}

We deduce a practical sufficient condition for the non-vanishing of cyclotomic multiple harmonic sums.

\begin{Corollary} Assume that $\mathbb{Q}(\xi_{i})$ is of degree $l_{i}\geq m - m_{0}$ with $l_{i}$ and $l'_{i}$ coprime for all $i\not=i'$.
Then, the multiple harmonic sum $\displaystyle\sum_{m_{0}<m_{1}<\ldots <m_{d}<m} \frac{\xi_{1}^{m_{1}} \cdots \xi_{d}^{m_{d}}}{m_{1}^{n_{1}} \cdots m_{d}^{n_{d}}}$ and all its subsums are non-zero.
\end{Corollary}

\begin{proof} It suffices to check that we can apply proposition \ref{prop general field}.

By fact \ref{fact linearly disjoint coprime}, the extensions $\mathbb{Q}(\xi_{1}),\ldots,\mathbb{Q}(\xi_{d})$ are linearly disjoint over $\mathbb{Q}$. Thus, by fact \ref{fact linearly disjoint property} (d) we deduce 
$[\mathbb{Q}(\xi_{1}\cdots \xi_{d}):\mathbb{Q}(\xi_{1}\cdots \widehat{\xi_{i}} \cdots \xi_{d})] = [\mathbb{Q}(\xi_{i}):\mathbb{Q}] =l_{i}\geq m - m_{0}$.
 
Then, by applying fact \ref{fact linearly disjoint property} (a) and (b), we deduce for all $1\leq i \leq l \leq j\leq d$, we have
$[\mathbb{Q}(\xi_{i}\cdots \xi_{j}):\mathbb{Q}(\xi_{i}\cdots \widehat{\xi_{l}} \cdots \xi_{j})] = [\mathbb{Q}(\xi_{1}\cdots \xi_{d}):\mathbb{Q}(\xi_{1}\cdots \widehat{\xi_{l}} \cdots \xi_{j})]$, and for all $m_{0}\leq m'_{0} < m' \leq m$, we have $m' - m'_{0} \geq m - m_{0}$. In particular, for all $1\leq i \leq l \leq j\leq d$ and $m_{0}\leq m'_{0} < m' \leq m$, we have $[\mathbb{Q}(\xi_{i}\cdots \xi_{j}):\mathbb{Q}(\xi_{i}\cdots \widehat{\xi_{l}} \cdots \xi_{j})] \geq m' - m'_{0}$. Thus the desired inequality for the subsums and we can apply proposition \ref{prop general field}.
\end{proof}

\begin{Remark} By fact 2.6 (d), the hypothesis of corollary 2.9 implies that $[\mathbb{Q}(\xi_{1}\cdots \xi_{d}):\mathbb{Q}]\geq (m-m_{0})^{d}$.
\end{Remark}

\subsection{Bounds on the degrees of $\xi_{1}$ and $\xi_{d}$}

Among all the $P_{i}$'s (in the sense of fact 2.1), the case of $P_{1}$ and $P_{d}$ are particular :

\begin{Proposition} \label{P1 and Pd}For any multiple harmonic sum $h_{m_{0},m}((n_{i})_{d};(\xi_{i})_{d})$ we always have $P_{1} \not= 0 $ and $ P_{d} \not= 0$. In particular, if $[\mathbb{Q}(\xi_{1}\cdots \xi_{d}):\mathbb{Q}(\xi_{2}\cdots \xi_{d})]\geq m - m_{0}$ or $[\mathbb{Q}(\xi_{1}\cdots\xi_{d}):\mathbb{Q}(\xi_{1}\cdots \xi_{d-1})] \geq m- m_{0}$, a multiple harmonic sum $h_{m_{0},m}((n_{i})_{d};(\xi_{i})_{d})$ is non-zero.
\end{Proposition}

\begin{proof} We make the proof for $P_{1}$, the proof for $P_{d}$ is similar.

Assume $P_{1}=0$. This implies for all $m_{1}$ such that $m_{0}<m_{1}<m$, 
$\displaystyle \sum_{m_{1}<m_{2}<\ldots <m_{d}<m} \frac{\xi_{2}^{m_{2}} \cdots \xi_{d}^{m_{d}}}{m_{1}^{n_{1}} \cdots m_{d}^{n_{d}}}=0$, thus
$\displaystyle \sum_{m_{1}<m_{2}<\ldots <m_{d}<m} \frac{\xi_{2}^{m_{2}} \cdots \xi_{d}^{m_{d}}}{m_{2}^{n_{2}} \cdots m_{d}^{n_{d}}}=0$.
In particular, for $m_{1}> m_{0}+1$, for all $m_{2}$ such that $m_{1}<m_{2}<m_{d}$,
$\displaystyle\bigg( \sum_{m_{1}<m_{2}<\ldots <m_{d}<m} - \sum_{m_{1}-1<m_{2}<\ldots <m_{d}<m} \bigg) \frac{\xi_{2}^{m_{2}} \cdots \xi_{d}^{m_{d}}}{m_{2}^{n_{2}} \cdots m_{d}^{n_{d}}} =  \sum_{m_{1}=m_{2}<m_{3}\ldots <m_{d}<m} \frac{\xi_{2}^{m_{2}} \cdots \xi_{d}^{m_{d}}}{m_{2}^{n_{2}} \cdots m_{d}^{n_{d}}}= 0$.

Thus, for all $m_{2}>m_{0}+2$,
$\displaystyle \sum_{m_{2}<m_{3} \ldots <m_{d}<m} \frac{\xi_{3}^{m_{3}} \cdots \xi_{d}^{m_{d}}}{m_{3}^{n_{3}} \cdots m_{d}^{n_{d}}} = 0 $.
Similarly, we show by induction on $i$ that for all $m_{i}>m_{0}+i$, 

$$ \sum_{m_{i}<m_{i+1}< \ldots <m_{d}<m} \frac{\xi_{i+1}^{m_{i+1}} \cdots \xi_{d}^{m_{d}}}{m_{i}^{n_{i}} \cdots m_{d}^{n_{d}}} = 0 . $$ 

Applying this for $i=d$ we find, for all $m_{d}> m_{0} + d$,
$ \frac{\xi_{d}^{m_{d}}}{m_{d}^{n_{d}}} = 0 $
whence a contradiction. 

We are left with the case where there does not exist $(m_{1},\ldots,m_{d})$ such that $m_{0}<m_{1}<\ldots<m_{d}<m$ and $m_{d}>m_{0}+d$. This is the case where $m = m_{0}+d+1$ and the only $(m_{1},\ldots,m_{d})$ such that $m_{0}<m_{1}<\ldots<m_{d}<m$ is $(m_{0}+i)_{i}$. In that case, $P_{1}$ is a non-zero multiple of a monomial, in particular it is non-zero.
\end{proof}

\subsection{The case of depth $\leq 2$}

In depth $\leq 2$, proposition \ref{P1 and Pd} allows to remove the hypothesis $P_{i}\not=0$ in fact \ref{fact 2.1} without cost :

\begin{Corollary} \label{corollary depth 2} (i) The harmonic sum $\displaystyle\sum_{m_{0}<m_{1}<m} \frac{\xi^{m_{1}}}{m_{1}^{n_{1}}}$ is non-zero if $[\mathbb{Q}(\xi):\mathbb{Q}] \geq m - m_{0}$. In particular, given $m_{0},m,n_{1}$, there are finitely many $\xi$'s such that $\displaystyle\sum_{m_{0}<m_{1}<m} \frac{\xi^{m_{1}}}{m_{1}^{n_{1}}}=0$.

(ii) The multiple harmonic sum
$\displaystyle\sum_{m_{0}<m_{1}<m_{2}<m} \frac{\xi_{1}^{m_{1}} \xi_{2}^{m_{2}}}{m_{1}^{n_{1}} m_{2}^{n_{2}}}$ is non-zero if  
$[\mathbb{Q}(\xi_{1},\xi_{2}):\mathbb{Q}(\xi_{2})] \geq m- m_{0}$ or 
$[\mathbb{Q}(\xi_{1},\xi_{2}):\mathbb{Q}(\xi_{1})] \geq m- m_{0}$. In particular, given $m_{0},m,n_{1},n_{2}$, there are finitely many $(\xi_{1},\xi_{2})$ such that $\displaystyle\sum_{m_{0}<m_{1}<m_{2}<m} \frac{\xi_{1}^{m_{1}} \xi_{2}^{m_{2}}}{m_{1}^{n_{1}} m_{2}^{n_{2}}}=0$
\end{Corollary}

\begin{proof} Immediate from proposition \ref{P1 and Pd} and fact \ref{fact 2.1}.
\end{proof}

Here is how (ii) above applies to different situations.

\begin{Example} If $\mathbb{Q}(\xi_{1}),\mathbb{Q}(\xi_{2})$ are linearly disjoint, we have, by fact 2.6 (d), 
$[\mathbb{Q}(\xi_{1},\xi_{2}):\mathbb{Q}] = [\mathbb{Q}(\xi_{1},\xi_{2}):\mathbb{Q}(\xi_{2})][\mathbb{Q}(\xi_{1},\xi_{2}):\mathbb{Q}(\xi_{1})] $. Thus if the multiple harmonic sum of corollary \ref{corollary depth 2} (ii) is zero we have 
$$ [\mathbb{Q}(\xi_{1},\xi_{2}):\mathbb{Q}] \leq (m-m_{0} - 1)^{2} . $$
\end{Example}

\begin{Example} If we have $\xi_{2} = \xi_{1}^{a}$, denoting $\xi_{1}= \xi$, and denoting by $N$ the order of $\xi$, we have  
$$ [\mathbb{Q}(\xi_{1},\xi_{2}):\mathbb{Q}(\xi_{1})] = 1, $$
$$ [\mathbb{Q}(\xi_{1},\xi_{2}):\mathbb{Q}(\xi_{2})]=[\mathbb{Q}(\xi):\mathbb{Q}(\xi^{a})] = \frac{\varphi(N)}{\varphi(\frac{N}{\gcd(a,N)})}. $$ 

(Indeed let us write $\xi = e^{2i\pi \frac{u}{N}}$ with $u$ and $N$ coprime, we have $\xi^{a} = e^{2i\pi \frac{ua}{N}}= e^{2i\pi \frac{ua'}{N'}}$, where $a = \gcd(a,N) a'$ and 
$N = \gcd(a,N)N'$, and we can deduce that the order of $\xi^{a}$ is $N' = \frac{N}{\gcd(a,N)}$, thus 
$[\mathbb{Q}(\xi^{a}):\mathbb{Q}]=\varphi(\frac{N}{\gcd(a,N)})$.) So if the multiple harmonic sum of corollary \ref{corollary depth 2} (ii) is zero, we have 
$$ \frac{\varphi(N)}{\varphi(\frac{N}{\gcd(a,N)})} \leq m- m_{0} - 1 . $$

Noting that for all integer $l$ we have $\sqrt{\frac{l}{2}} \leq \varphi(l) \leq l$, whence $\frac{\varphi(N)}{\varphi(\frac{N}{\gcd(a,N)})} \geq \frac{\gcd(a,N)}{\sqrt{2N}}$, thus the inequality implies a bound on $\gcd(a,N)$ : 

$$ \gcd(a,N) \leq (m-m_{0}-1) \sqrt{2N} . $$ 

If $a |N$, and $N = ab$ with $a$ and $b$ coprime, then we have $\displaystyle \frac{\varphi(N)}{\varphi(\frac{N}{\gcd(a,N)})} = \frac{\varphi(N)}{\varphi(b)} = \varphi(a)$ by multiplicativity of $\varphi$. Thus the inequality is $\varphi(a) \leq m-1$. It implies a bound on a : 
$$ \sqrt{\frac{a}{2}} \leq m-1 . $$
\end{Example}

\begin{Example} In general we have $\xi_{1} = \xi^{a_{1}}$, $\xi_{2} = \xi^{a_{2}}$ for a root of unity $\xi$ and $a_{1},a_{2} \in \mathbb{N}$. 

Then we have $\mathbb{Q}(\xi_{1},\xi_{2}) = \mathbb{Q}(\xi^{\gcd(a_{1},a_{2})})$. Let $\xi'=\xi^{\gcd(a_{1},a_{2})}$, and 
$a_{i}'=\frac{a_{i}}{\gcd(a_{1},a_{2})}$, we have $\xi_{i} = {\xi'}^{a_{i}'}$. Let $N'$ the order of $\xi'$. We have $N' = \frac{N}{\gcd(a1,a2,N)}$, $[\mathbb{Q}(\xi_{1},\xi_{2}) : \mathbb{Q}(\xi_{2})] = [\mathbb{Q}({\xi'} : \mathbb{Q}({\xi'}^{\frac{a_{1}}{\gcd(a_{1},a_{2})}})] = \frac{\varphi(N')}{\varphi(\frac{N'}{\gcd(a_{1}',N')})}$, $[\mathbb{Q}(\xi_{1},\xi_{2}) : \mathbb{Q}(\xi_{1})] =  \frac{\varphi(N')}{\varphi(\frac{N'}{\gcd(a_{2}',N')})}$.
So if the multiple harmonic sum of corollary \ref{corollary depth 2} (ii) is zero, we have 
$$\forall i=1,2,  \frac{\varphi(N')}{\varphi(\frac{N'}{\gcd(a_{i}',N')})} \leq m - m_{0} - 1 . $$
\end{Example}

\subsection{Use of the distribution relation}

In a different direction, we prove that the non-vanishing of a multiple harmonic sum implies the non-vanishing of other multiple harmonic sums.

\begin{Proposition} Assume that a multiple harmonic sum $\tilde{h}_{m}((n_{i})_{d};(\xi_{i})_{d+1})$ is non-zero. Let $M$ be a positive integer. For any $1\leq i\leq d+1$ let $\tilde{\xi}_{i}$ be a $M$-th root of $\xi_{i}$. Then, there exists $M$-th roots of unity $\rho_{i}$ ($1 \leq i \leq d+1$) such that 
$$ \tilde{h}_{Mm} \big( (n_{i})_{d};(\rho_{i}\tilde{\xi}_{i})_{d+1} \big) \not=0 . $$
\end{Proposition}

\begin{proof} This follows from the distribution relation for cyclotomic multiple harmonic sums : for any $d \in \mathbb{N}^{\ast}$, any positive integers $n_{i}$ ($1 \leqslant i \leqslant d$) and roots of unity $\xi_{i}$ ($1 \leqslant i \leqslant d+1$), and for $M$ an integer, we have
$$ M^{\sum_{i=1}^{d}(n_{i}-1)} \displaystyle\sum_{\rho_{1}^{M}=1,\ldots,\rho_{d+1}^{M}=1} \tilde{h}_{m'} \big( (n_{i})_{d};(\rho_{i}\xi_{i})_{d+1} \big) = \left\{ \begin{array}{ll} \tilde{h}_{\frac{m'}{M}}\big((n_{i})_{d};(\xi_{i}^{M})_{d+1}\big) & \text{ }\text{if} \text{ } M\text{ }|\text{ }m'
\\ 0 & \text{ }\text{if}\text{ } M \nmid m'
\end{array}
\right. . $$

This is a standard fact but let us review its proof. Let $N$ be an integer divisible by $M$ such that the $\xi_{i}$'s are $N$-th roots of unity. For a sequence of global differential forms $\omega_{1},\ldots,\omega_{n}$ on $\mathbb{P}^{1} - \{0,\mu_{N},\infty\}$, such that $\omega_{1}$ has no pole at $0$, the formal iterated integral $\I(\omega_{1},\ldots,\omega_{n}) \in K[[z]]$ is defined by induction by 
$\I(\emptyset)=1$ and, for $n\geqslant 1$, $d\I(\omega_{1},\ldots,\omega_{n}) = \I(\omega_{1},\ldots,\omega_{n-1})\omega_{n}$ and $\I(\omega_{1},\ldots,\omega_{n})(0)=0$. Let us consider a sequence of differential forms as follows, where the $\xi_{i}$'s are $N$-th roots of unity
$\displaystyle (\omega_{1},\ldots,\omega_{n}) = \bigg(\frac{d(z^{M})}{z^{M}-\xi_{1}^{M}},\underbrace{\frac{d(z^{M})}{z^{M}},\ldots,\frac{d(z^{M})}{z^{M}}}_{n_{1}-1},\ldots,
\frac{d(z^{M})}{z^{M}-\xi_{d}^{M}},\underbrace{\frac{d(z^{M})}{z^{M}},\ldots,\frac{d(z^{M})}{z^{M}}}_{n_{d}-1},\frac{d (z^{M})}{z^{M}-\xi_{d+1}^{M}}\displaystyle\bigg)$. We have $\I(\omega_{1},\ldots,\omega_{n})=\displaystyle\sum_{0<m_{1}<\ldots<m_{d}<m} \frac{(\frac{\xi_{2}^{M}}{\xi_{1}^{M}})^{ m_{1}} \ldots (\frac{\xi_{d+1}^{M}}{\xi^{M}_{d}})^{m_{d}} (\frac{z^{M}}{\xi^{M}_{d+1}})^{m}}{m_{1}^{n_{1}}\ldots m_{d}^{n_{d}}m}$. On the other hand, we have $\displaystyle \frac{d(z^{M})}{z^{M}}=M\frac{dz}{z}$ and
$\displaystyle \frac{d(z^{M})}{z^{M}-\xi^{M}} = \frac{d(z^{M}\xi^{-M})}{z^{M}\xi^{-M}-1}=\sum_{\rho^{M}=1} \frac{d(z\xi^{-1})}{z\xi^{-1} - \rho} = \sum_{\rho^{N}=1} \frac{dz}{z - \rho \xi}$, thus 
\newline 
$\I(\omega_{1},\ldots,\omega_{n}) = M^{\sum_{i=1}^{d}(n_{i}-1)} \displaystyle\sum_{\rho_{1}^{M}=1,\ldots,\rho_{d+1}^{M}=1}\sum_{0<m_{1}<\ldots<m_{d}<m} \frac{(\frac{\rho_{2}\xi_{2}}{\rho_{1}\xi_{1}})^{m_{1}} \ldots (\frac{\rho_{d+1}\xi_{d+1}}{\rho_{d}\xi_{d}})^{m_{d}}(\frac{z}{\rho_{d+1}\xi_{d+1}})^{m}}{m_{1}^{n_{1}}\ldots m_{d}^{n_{d}}m}$. The equality between the two expressions of $\I(\omega_{1},\ldots,\omega_{n})$ gives the result.

Now for any $\tilde{h}_{m}((n_{i})_{d};(\xi_{i})_{d})\not=0$, and for any positive integer $M$, denoting by $m'=Mm$, we can find $\tilde{\xi}_{i}$ such that $\xi_{i} = \tilde{\xi}_{i}^{M}$. By the distribution relation we have

$$ M^{\sum_{i=1}^{d}(n_{i}-1)} \displaystyle\sum_{\rho_{1}^{M}=1,\ldots,\rho_{d+1}^{M}=1} h_{Mm} \big( (n_{i})_{d};(\rho_{i}\tilde{\xi}_{i})_{d+1} \big) = h_{m} \big((n_{i})_{d};(\xi_{i})_{d+1}\big) . $$
\end{proof}

\section{Conditions on $(n_{i})_{d}$}

We now find sufficients conditions on $(n_{i})_{d}$ to have 
$h_{m}((n_{i})_{d};(\xi_{i})_{d})\not=0$. The idea is a complex analogue of the $p$-adic idea in \S1 : we look at the complex absolute value of a cyclotomic multiple harmonic sum, and we see that, in certain cases, one term of the sum has a bigger absolute value than the sum of the other terms. This implies that the cyclotomic multiple harmonic sum is non-zero.

\subsection{For multiple harmonic sums}

In this section we prove that, given $m$ and $d$, there are only finitely many $n$'s such that there exists $(n_{i})_{d}$, $(\xi_{i})_{d}$ with $n_{d}=n$, and $h_{m}((n_{i})_{d};(\xi_{i})_{d}) = 0$, i.e. we prove that for $n_{d}$ large enough, with an explicit bound, we always have $h_{m}((n_{i})_{d};(\xi_{i})_{d}) \not= 0$. Then we prove an analogous but conditional result on $n_{i}$ with $i<d$.

\begin{Proposition} \label{prop nd}Assume that $\displaystyle n_{d}> \frac{\log\big({m-m_{0}-1 \choose d} - 1\big)}{ \log(\frac{m_{0}+d+1}{m_{0}+d})}$. Then for any positive integers $n_{i}$ ($1\leqslant i \leqslant d-1$), and for any roots of unity $\xi_{i}$ ($1\leqslant i \leqslant d+1$), we have  
$h_{m_{0},m} \big((n_{i})_{d};(\xi_{i})_{d}\big) \not=0$.
\end{Proposition}

\begin{proof} Denoting by $S_{m,d}$ the set of tuples $(m_{1},\ldots,m_{d}) \in \mathbb{N}^{d}$ such that $m_{0}<m_{1}<\ldots<m_{d}<m$ and $(m_{1},\ldots,m_{d})\not=(m_{0}+1,\ldots,m_{0}+d)$, we have 
$$ h_{m} \big((n_{i})_{d};(\xi_{i})_{d} \big) = \frac{\xi_{1}^{m_{0}+1} \ldots \xi_{d}^{m_{0}+d}}{(m_{0}+1)^{n_{1}} \ldots (m_{0}+d)^{n_{d}}} + \sum_{(m_{1},\ldots,m_{d}) \in S_{m,d}}
\frac{\xi_{1}^{m_{1}} \ldots \xi_{d}^{m_{d}}}{m_{1}^{n_{1}} \ldots m_{d}^{n_{d}}} . $$
For any $(m_{1},\ldots,m_{d}) \in S_{m,d}$ we have $m_{d}\geqslant m_{0}+d+1$ and 
$m_{1}^{n_{1}} \ldots m_{d}^{n_{d}} \geqslant (m_{0}+d+1)^{n_{d}} \prod_{i=1}^{d-1}(m_{0}+i)^{n_{i}}$, whence 
$\displaystyle \bigg| \sum_{(m_{1},\ldots,m_{d}) \in S_{m,d}}
\frac{\xi_{1}^{m_{1}} \ldots \xi_{d}^{m_{d}}}{m_{1}^{n_{1}} \ldots m_{d}^{n_{d}}} \bigg| \leqslant  
\frac{|S_{m,d}|}{(m_{0}+d+1)^{n_{d}} \prod_{i=1}^{d-1}(m_{0}+i)^{n_{i}}}$. On the other hand, we have 
$\displaystyle \bigg| \frac{\prod_{i=1}^{d}\xi_{i}^{m_{0}+i}  }{\prod_{i=1}^{d}(m_{0}+i)^{n_{i}}} \bigg| = \frac{1}{\prod_{i=1}^{d}(m_{0}+i)^{n_{i}}}$. Finally, we have the equivalences 
$\displaystyle\frac{|S_{m,d}|}{(m_{0}+d+1)^{n_{d}}\prod_{i=1}^{d-1}(m_{0}+i)^{n_{i}}} < \frac{1}{\prod_{i=1}^{d}(m_{0}+i)^{n_{i}}} \Leftrightarrow |S_{m,d}| < \bigg(\frac{m_{0}+d+1}{m_{0}+d} \bigg)^{n_{d}} \Leftrightarrow \frac{\log(|S_{m,d}|)}{\log(\frac{m_{0}+d+1}{m_{0}+d})} < n_{d}$, and $\displaystyle |S_{m,d}| = {m-m_{0}-1 \choose d} - 1$.
\end{proof}

\begin{Proposition} \label{prop ni}Let $1 \leq i \leq d$. Let us assume that $\displaystyle \sum\limits_{m_{0}+i+2\leq m_{i+1} < \cdots < m_{d} < m} \frac{\xi_{i+1}^{m_{i+1}} \cdots \xi_{d}^{m_{d}}}{m_{i+1}^{n_{i+1}} \cdots m_{d}^{n_{d}}}$ and
\newline $\displaystyle \sum\limits_{m_{0}+i+1\leq  m_{i+1} < \cdots < m_{d} < m} \frac{\xi_{i+1}^{m_{i+1}} \cdots \xi_{d}^{m_{d}}}{m_{i+1}^{n_{i+1}} \cdots m_{d}^{n_{d}}}$ are non-zero.
Then, for $n_{i}> \frac{\log  \Bigg( \displaystyle\frac{\bigg| \sum\limits_{m_{0}+i+2\leq m_{i+1} < \cdots < m_{d} < m} \frac{\xi_{i+1}^{m_{i+1}} \cdots \xi_{d}^{m_{d}}}{m_{i+1}^{n_{i+1}} \cdots m_{d}^{n_{d}}} \bigg|}{\bigg| \sum\limits_{m_{0}+i+1\leq  m_{i+1} < \cdots < m_{d} < m} \frac{\xi_{i+1}^{m_{i+1}} \cdots \xi_{d}^{m_{d}}}{m_{i+1}^{n_{i+1}} \cdots m_{d}^{n_{d}}} \bigg|} \Bigg)}{\log\big( \frac{i+1}{i}\big)}$, $h_{m_{0},m}((n_{i})_{d},(\xi_{i})_{d})$ is non-zero.
\end{Proposition}

\begin{proof} We write $h_{m}((n_{i})_{d},(\xi_{i})_{d})$ as the sum of two terms corresponding to the two subdomains $m_{i}=m_{0}+i$ and $m_{i} > m_{0}+i$ of the domain of summation.
	
The $m_{i}=m_{0}+i$ term is $\displaystyle\frac{\xi_{1}^{m_{0}+1}\cdots \xi_{i}^{m_{0}+i}}{(m_{0}+1)^{n_{1}} \cdots (m_{0}+i)^{n_{i}}} \sum_{m_{0}+i+1\leq m_{i+1} < \cdots < m_{d} < m} \frac{\xi_{i+1}^{m_{i+1}} \cdots \xi_{d}^{m_{d}}}{m_{i+1}^{n_{i+1}} \cdots m_{d}^{n_{d}}}$.

The $m_{i}>m_{0}+i$ term is $\displaystyle\frac{\xi_{1}^{m_{0}+1}\cdots \xi_{i-1}^{m_{0}+i-1}}{(m_{0}+1)^{n_{1}} \cdots (m_{0}+i-1)^{n_{i-1}}} \sum_{m_{0}+i+1\leq m_{i}< m_{i+1} < \cdots < m_{d} < m} \frac{\xi_{i}^{m_{i}} \cdots \xi_{d}^{m_{d}}}{m_{i}^{n_{i}} \cdots m_{d}^{n_{d}}}$.

If we have the inequality |$m_{i}>m_{0}+i$ term| < |$m_{i}=m_{0}+i$ term|, then we have $h_{m_{0},m}((n_{i})_{d},(\xi_{i})_{d})$.

We have
\begin{multline*}\displaystyle \bigg|\frac{\xi_{1}^{m_{0}+1}\cdots \xi_{i-1}^{m_{0}+i-1}}{(m_{0}+1)^{n_{1}} \cdots (m_{0}+i-1)^{n_{i-1}}} \sum_{m_{0}+i+1\leq m_{i}< m_{i+1} < \cdots < m_{d} < m} \frac{\xi_{i}^{m_{i}} \cdots \xi_{d}^{m_{d}}}{m_{i}^{n_{i}} \cdots m_{d}^{n_{d}}}\bigg|
\\ \leq \frac{1}{(m_{0}+1)^{n_{1}} \cdots (m_{0}+i-1)^{n_{i-1}}} 
\frac{1}{(m_{0}+i+1)^{n_{i}}}
\bigg| \sum_{i+2\leq m_{i+1} < \cdots < m_{d} < m} \frac{\xi_{i+1}^{m_{i+1}} \cdots \xi_{d}^{m_{d}}}{m_{i+1}^{n_{i+1}} \cdots m_{d}^{n_{d}}} \bigg|.
\end{multline*}

Thus, the inequality |$m_{i}>m_{0}+i$ term| < |$m_{i}=m_{0}+i$ term| is implied by 

\begin{multline*}\frac{1}{(m_{0}+1)^{n_{1}} \cdots (m_{0}+i-1)^{n_{i-1}}} 
\frac{1}{(m_{0}+i+1)^{n_{i}}}|\sum_{m_{0}+i+2\leq m_{i+1} < \cdots < m_{d} < m} \frac{1}{m_{i+1}^{n_{i+1}} \cdots m_{d}^{n_{d}}}|
\\ < \frac{1}{(m_{0}+1)^{n_{1}} \cdots (m_{0}+i)^{n_{i}}} \bigg| \sum_{m_{0}+i+1\leq  m_{i+1} < \cdots < m_{d} < m} \frac{\xi_{i+1}^{m_{i+1}} \cdots \xi_{d}^{m_{d}}}{m_{i+1}^{n_{i+1}} \cdots m_{d}^{n_{d}}} \bigg|,
\end{multline*}

i.e. 

$$ \bigg( \frac{m_{0}+i}{m_{0}+i+1}\bigg)^{n_{i}}
\bigg| \sum_{m_{0}+i+2\leq m_{i+1} < \cdots < m_{d} < m} \frac{\xi_{i+1}^{m_{i+1}} \cdots \xi_{d}^{m_{d}}}{m_{i+1}^{n_{i+1}} \cdots m_{d}^{n_{d}}} \bigg| < \bigg| \sum_{m_{0}+i+1\leq  m_{i+1} < \cdots < m_{d} < m} \frac{\xi_{i+1}^{m_{i+1}} \cdots \xi_{d}^{m_{d}}}{m_{i+1}^{n_{i+1}} \cdots m_{d}^{n_{d}}} \bigg| . $$

Assuming the multiple harmonic sum in the left-hand side above is non-zero, this is equivalent to 

\begin{equation} \label{eq:this inequality} \bigg( \frac{m_{0}+i}{m_{0}+i+1}\bigg)^{n_{i}}
< \displaystyle\frac{\bigg|\sum\limits_{m_{0}+i+1\leq  m_{i+1} < \cdots < m_{d} < m} \frac{\xi_{i+1}^{m_{i+1}} \cdots \xi_{d}^{m_{d}}}{m_{i+1}^{n_{i+1}} \cdots m_{d}^{n_{d}}} \bigg|}{\bigg| \sum\limits_{m_{0}+i+2\leq m_{i+1} < \cdots < m_{d} < m} \frac{\xi_{i+1}^{m_{i+1}} \cdots \xi_{d}^{m_{d}}}{m_{i+1}^{n_{i+1}} \cdots m_{d}^{n_{d}}} \bigg|} .
\end{equation}
 
If $\displaystyle\frac{\bigg| \sum\limits_{m_{0}+i+1\leq  m_{i+1} < \cdots < m_{d} < m} \frac{\xi_{i+1}^{m_{i+1}} \cdots \xi_{d}^{m_{d}}}{m_{i+1}^{n_{i+1}} \cdots m_{d}^{n_{d}}} \bigg|}{\bigg| \sum\limits_{m_{0}+i+2\leq m_{i+1} < \cdots < m_{d} < m} \frac{\xi_{i+1}^{m_{i+1}} \cdots \xi_{d}^{m_{d}}}{m_{i+1}^{n_{i+1}} \cdots m_{d}^{n_{d}}} \bigg|} \geq 1$, equation (\ref{eq:this inequality}) is true for all $n_{i}$, since $\big( \frac{i}{i+1}\big)^{n_{i}} \leq 1$.

Otherwise, equation (\ref{eq:this inequality}) is equivalent to 
$ n_{i} \log\big( \frac{m_{0}+i+1}{m_{0}+i}\big)
> \log  \Bigg( \displaystyle\frac{\bigg| \sum\limits_{m_{0}+i+2\leq m_{i+1} < \cdots < m_{d} < m} \frac{\xi_{i+1}^{m_{i+1}} \cdots \xi_{d}^{m_{d}}}{m_{i+1}^{n_{i+1}} \cdots m_{d}^{n_{d}}} \bigg|}{\bigg| \sum\limits_{m_{0}+i+1\leq  m_{i+1} < \cdots < m_{d} < m} \frac{\xi_{i+1}^{m_{i+1}} \cdots \xi_{d}^{m_{d}}}{m_{i+1}^{n_{i+1}} \cdots m_{d}^{n_{d}}} \bigg|} \Bigg) $

i.e. $n_{i} > \frac{\log  \Bigg( \displaystyle\frac{\bigg| \sum\limits_{m_{0}+i+2\leq m_{i+1} < \cdots < m_{d} < m} \frac{\xi_{i+1}^{m_{i+1}} \cdots \xi_{d}^{m_{d}}}{m_{i+1}^{n_{i+1}} \cdots m_{d}^{n_{d}}} \bigg|}{\bigg| \sum\limits_{m_{0}+i+1\leq  m_{i+1} < \cdots < m_{d} < m} \frac{\xi_{i+1}^{m_{i+1}} \cdots \xi_{d}^{m_{d}}}{m_{i+1}^{n_{i+1}} \cdots m_{d}^{n_{d}}} \bigg|} \Bigg)}{\log\big( \frac{m_{0}+i+1}{m_{0}+i}\big)}$.
\end{proof}

\begin{Remark} It is possible to combine proposition 3.2 with results of non-vanishing in depth 1, which do not require conditions on $n_{d}$, such as the results of \S2.4, and obtain unconditional non-vanishing results with conditions on $n_{i}$ with $i<d$, by decreasing induction on $i$.
\end{Remark}

\subsection{For cyclotomic multiple zeta values}

Actually, the above methods can also be used to prove the non-vanishing of a CMZV.

\begin{Proposition} Let $(n_{i})_{d}$ such that $n_{d}>\frac{\log \bigg(d + \zeta(n_{1},\ldots,n_{d-1}) (d+1)\prod_{i=1}^{d-1} i^{n_{i}}\bigg)}{\log \big(\frac{d+1}{d}\big)}$, then for all $(\xi'_{i})_{d}$ we have $\zeta((n_{i})_{d},(\xi'_{i})_{d})\not=0$.
\end{Proposition}

\begin{proof} For any $m$, we have 
$\displaystyle \zeta((n_{i})_{d},(\xi'_{i})_{d}) = \sum_{0<m_{1}<\cdots<m_{d}} \frac{\xi_{1}^{m_{1}} \cdots \xi_{d}^{m_{d}}}{m_{1}^{n_{1}} \cdots m_{d}^{n_{d}}} 
= h_{m}((n_{i})_{d};(\xi_{i})_{d}) + \sum_{\substack{0<m_{1}<\cdots<m_{d} \\ m_{d}\geq m}} \frac{\xi_{1}^{m_{1}} \cdots \xi_{d}^{m_{d}}}{m_{1}^{n_{1}} \cdots m_{d}^{n_{d}}}$, with $\xi_{i} = \frac{\xi'_{i+1}}{\xi'_{i}}$, where $\xi'_{d+1}=1$. By the proof of proposition 3.1, we write

$$ \zeta((n_{i})_{d},(\xi_{i})_{d}) = \frac{\xi_{1}^{1} \cdots \xi_{d}^{d}}{m_{1}^{n_{1}} \cdots m_{d}^{n_{d}}} + 
\sum_{\substack{0<m_{1}<\cdots<m_{d}<m \\ (m_{i})_{d} \in  S_{m,d}}} \frac{\xi_{1}^{m_{1}} \cdots \xi_{d}^{m_{d}}}{m_{1}^{n_{1}} \cdots m_{d}^{n_{d}}} + \sum_{\substack{0<m_{1}<\cdots<m_{d} \\ m_{d}\geq m}} \frac{\xi_{1}^{m_{1}} \cdots \xi_{d}^{m_{d}}}{m_{1}^{n_{1}} \cdots m_{d}^{n_{d}}} $$

We bound the two last terms above: $$
\begin{array}{ll}
\displaystyle \bigg|\sum_{\substack{0<m_{1}<\cdots<m_{d} \\ m_{d}\geq m}} \frac{\xi_{1}^{m_{1}} \cdots \xi_{d}^{m_{d}}}{m_{1}^{n_{1}} \cdots m_{d}^{n_{d}}}\bigg|
& \displaystyle \leq \zeta(n_{1},\ldots,n_{d-1}) \sum_{m_{d}\geq m} \frac{1}{m_{d}^{n_{d}}} 
\\ & \displaystyle \leq \zeta(n_{1},\ldots,n_{d-1}) \int_{m-1}^{\infty} \frac{dt}{t^{n_{d}}} = \zeta(n_{1},\ldots,n_{d-1})\frac{1}{(n_{d}-1) (m-1)^{n_{d}-1}}
\end{array}, $$
and 
$$ \bigg|\sum_{\substack{0<m_{1}<\cdots<m_{d}<m \\ (m_{i})_{d} \in  S_{m,d}}} \frac{\xi_{1}^{m_{1}} \cdots \xi_{d}^{m_{d}}}{m_{1}^{n_{1}} \cdots m_{d}^{n_{d}}} \bigg| \leq \frac{{m-1 \choose d} - 1}{(d+1)^{n_{d}} \prod_{i=1}^{d-1} i^{n_{i}}} . $$

Thus, the non-vanishing of is implied by the inequality 

$$ \frac{{m-1 \choose d} - 1}{(d+1)^{n_{d}} \prod_{i=1}^{d-1} i^{n_{i}}} + \zeta(n_{1},\ldots,n_{d-1})\frac{1}{(n_{d}-1) (m-1)^{n_{d}-1}} < \frac{1}{\prod_{i=1}^{d} i^{n_{i}}}, $$

which is equivalent to 

$$ \bigg({m-1  \choose d} - 1\bigg) \big( \frac{d}{d+1}\big)^{n_{d}} 
+ \zeta(n_{1},\ldots,n_{d-1}) \frac{(\prod_{i=1}^{d-1} i^{n_{i}}) d^{n_{d}}}{(n_{d}-1) (m-1)^{n_{d}-1}} < 1, $$

i.e.

$$ \bigg({m-1  \choose d} - 1\bigg) \big( \frac{d}{d+1}\big)^{n_{d}} 
+ \zeta(n_{1},\ldots,n_{d-1}) \frac{(\prod_{i=1}^{d-1} i^{n_{i}})(m-1)}{(n_{d}-1)} \big(\frac{d}{m-1}\big)^{n_{d}}  < 1 . $$

This is implied by 

$$ \max(\big(\frac{d}{m-1}\big)^{n_{d}},\big( \frac{d}{d+1}\big)^{n_{d}}) < \frac{1}{({m-1  \choose d} - 1) + \zeta(n_{1},\ldots,n_{d-1}) \frac{(m-1)\prod_{i=1}^{d-1} i^{n_{i}}}{n_{d}-1}} . $$

Necessarily we must take $m>d+1$ thus $m-1 \geq d+1$. Thus the inequality amounts to 
$$ \big( \frac{d}{d+1}\big)^{n_{d}} < ({m-1  \choose d} - 1) + \zeta(n_{1},\ldots,n_{d-1}) \frac{(m-1)\prod_{i=1}^{d-1} i^{n_{i}}}{(n_{d}-1)} ,$$
i.e. 
$$ n_{d} > \frac{\log \bigg( \big({m-1 \choose d} - 1\big) + \zeta(n_{1},\ldots,n_{d-1}) \frac{(m-1)\prod_{i=1}^{d-1} i^{n_{i}}}{n_{d}-1} \bigg)}{\log \big(\frac{d+1}{d}\big)} . $$

The right-hand side is an increasing function of $m\geq d+2$. Thus, to have the optimal bound, we choose $m=d+2$. The inequality that we want is thus 
$$ n_{d} > \frac{\log \bigg( d + \zeta(n_{1},\ldots,n_{d-1}) \frac{(d+1)\prod_{i=1}^{d-1} i^{n_{i}}}{n_{d}-1} \bigg)}{\log \big(\frac{d+1}{d}\big)} . $$
The right-hand side depends on $n_{d}$ but is a decreasing function of $n_{d}$, and we have
$$ \frac{\log \bigg( d + \zeta(n_{1},\ldots,n_{d-1}) \frac{(d+1)\prod_{i=1}^{d-1} i^{n_{i}}}{n_{d}-1} \bigg)}{\log \big(\frac{d+1}{d}\big)} < \frac{\log \bigg( d + \zeta(n_{1},\ldots,n_{d-1}) (d+1)\prod_{i=1}^{d-1} i^{n_{i}} \bigg)}{\log \big(\frac{d+1}{d}\big)} . $$
Whence the result.
\end{proof}

\begin{Remark} It is also possible to write an analogue of proposition 3.2 for cyclotomic multiple zeta values.
\end{Remark}

\section{The case of $N=2$}

Let a multiple harmonic sum $h_{m}((n_{i})_{d},(\xi_{i})_{d})$ such that $\xi_{i} = \pm 1$ for all $i$. In particular, $h_{m}((n_{i})_{d},(\xi_{i})_{d})$ is in $\mathbb{R}$ and is expressed as an iterated alternating sum as follows : let $i_{1},\ldots,i_{r}$ be the integers such that $\xi_{i_{1}} = \cdots \xi_{i_{r}} = -1$, with $1 \leq r \leq d$ and $1 \leq i_{1}<\ldots < i_{r} \leq d$ ; for $M> 0$, let

$$ u(M) = \sum_{\substack{m_{i_{1}}+\cdots+m_{i_{r}}=M \\ 0<m_{1}<\ldots<m_{d}<m}} \frac{1}{m_{1}^{n_{1}} \cdots m_{d}^{n_{d}}} . $$

We have (the sum is finite) :

\begin{equation} \label{eq:zero} h_{m}((n_{i})_{d},(\xi_{i})_{d}) = \sum_{0<M} (-1)^{M}u(M) .
\end{equation}

\begin{Lemma} \label{lemma alternative}
The set $\{M\text{ }|\text{ }u(M) \not=0\}$ is the segment $[M_{\min}=\sum\limits_{j=1}^{r}i_{j},M_{\max}=r(m-d-1) + \sum\limits_{j=1}^{r} i_{j}]$ of $\mathbb{N}$. 
\end{Lemma}

\begin{proof} $u(M)$ is non-zero if and only if the sum $\sum\limits_{\substack{m_{i_{1}}+\cdots+m_{i_{r}}=M \\ 0<m_{1}<\ldots<m_{d}<m}} \frac{1}{m_{1}^{n_{1}} \cdots m_{d}^{n_{d}}}$ has a non-empty domain of summation, i.e. if and only if the set
$E= \{ (m_{1},\ldots,m_{d}) \in \mathbb{N}^{d} \text{ }|\text{ }m_{i_{1}}+\cdots+m_{i_{r}}=M,\text{ and }0<m_{1}<\ldots<m_{d}<m\}$ is non-empty.

The minimal element of $E$ is $M_{\min} = i_{1}+\cdots + i_{r}$, which corresponds to $m_{i}=i$ for all $i$.

The maximal element of $E$ is $M_{\max} =  \sum_{j=1}^{r} (m - (d-i_{j}) - 1) = r(m-d-1) + \sum_{j=1}^{r} i_{j}$. It corresponds to $m_{d-i} = m - i - 1$ for all $i$, i.e. $m_{i} = m - (d-i) - 1$ for all $i$. 

Let $M$ such that $u(M)$ is non-zero. Let $(m_{1},\ldots,m_{d})$ such that $0<m_{1}<\cdots<m_{d}<m$. Then we always have, for all $i$,
\begin{equation} \label{eq:theinequality} i \leq m_{i} \leq m - (d-i) - 1 .
\end{equation}
	
Now let $M$ be in $[M_{\min},M_{\max}]$. Let us find $(m_{1},\ldots,m_{d})$ in $E$ such that $m_{i_{1}}+ \cdots + m_{i_{r}}= M$, by induction on $M$.

Assume that $u(M-1)\not=0$. Let $(m_{1},\ldots,m_{d})$ in $E$ such that $m_{i_{1}}+ \cdots + m_{i_{r}}= M-1$.

We can assume that $M < M_{\max}$ because the property is already proved for $M_{\max}$.

Since $M < M_{\max}$, at least one of the inequalities $m_{i_{j}} \leq m - (d-i_{j}) - 1$ from equation (\ref{eq:theinequality}) is strict. We define $j_{0}= \max \{j\text{ }|\text{ }m_{i_{j}}< m - (d-i_{j}) - 1\}$. Let $j_{1}$ be the maximal $j$ such that for all $j_{0} \leq k \leq j_{1}$, we have $m_{j_{0}}+k=m_{j_{0}+k}$. Thus we have $m_{j_{1}+1} \geq m_{j_{1}}+2$. (If $j_{1}=d$, we denote by $m_{d+1}=m$.)

Consider the variant $(m'_{1},\ldots,m'_{d})$ of $(m_{1},\ldots,m_{d})$ obtained by replacing $m_{k}$ by $m_{k}+1$ for all $k \in [j_{0},j_{1}]$. 
Since $m_{j_{1}+1} \geq m_{j_{1}}+2$, we have
$m'_{j_{1}+1} = m'_{j_{1}}+1 < m_{j_{1}+1}=m'_{j_{1}+1}$, and  $0<m'_{1}<\ldots<m'_{d}<m$.

Moreover, we also have necessarily $i_{j_{0}+1}>j_{1}$: by definition of $j_{0}$, be have $m_{i_{j_{0}+1}} = m - (d-i_{j})-1$; this forces $m_{l} = m - (d-l)-1$ for all $l\geq i_{j_{0}+1}$; in particular, for all $l\geq i_{j_{0}+1}$, $l\leq d$, we have $m_{l+1} = m_{l}+1$. Whereas we have $m_{j_{1}+1} \geq m_{j_{1}}+2$.

Since $i_{j_{0}+1} > j_{1}$, the only $j$ such that $m'_{i_{j}} \not= m_{i_{j}}$ is $i_{j} = i_{j_{0}}$. And we have seen that, by definition, $m'_{i_{j_{0}}} = m_{i_{j_{0}}}+1$.

In particular, $m'_{i_{1}} + \cdots + m'_{i_{r}} = m_{i_{1}} + \cdots + m_{i_{r}} + 1=M-1+1=M$. 

As a conclusion, $(m'_{1},\ldots,m'_{d})$ is an element of $E$ such that $m'_{i_{1}} + \cdots + m'_{i_{r}}=M$. So we have $u(M)\not=0$.
\end{proof}

\begin{Lemma} \label{second lemma alternative}If $u$ is monotonous, then $u$ is decreasing. 
\end{Lemma}

\begin{proof} We have $\displaystyle u(M_{\min}) = \frac{1}{1^{m_{1}} \cdots d^{m_{d}}}$ and $\displaystyle u(M_{\max}) = \frac{1}{(m-d)^{m_{1}} \cdots (m-1)^{m_{d}}}$.
Since $m\geq d+1$, for all $1 \leq i \leq d$, we have $i \leq m-(d-i)-1 = (m-d-1) +i$, thus $i^{m_{i}} \leq (m-(d-i)-1)^{m_{i}}$, which implies $u(M_{\min}) \geq u(M_{\max})$.
\end{proof}

We note that $u$ is non-constant unless $m=d+1$ and $M_{\min}=M_{\max}$ : indeed, if $m>d+1$, the inequalities in the proof of lemma 4.2 are strict.

\begin{Proposition} \label{prop alternative}
Assume that $u$ is strictly monotonous on $[M_{\min},M_{\max}]$. Then, $h_{m}((n_{i})_{d},(\xi_{i})_{d})$ is non-zero.
\end{Proposition}

\begin{proof} We express $h_{m}((n_{i})_{d},(\xi_{i})_{d})$ by equation (\ref{eq:zero}) and lemma \ref{lemma alternative} : there are two cases :
\newline 
- first case: $M_{\max} - M_{\min}$ is odd:
$$ h_{m}((n_{i})_{d},(\xi_{i})_{d}) = 
(-1)^{M_{\min}} \big[ (u(M_{\min}) - u(M_{\min+1})) + \cdots + (u(M_{\max}-1) - u(M_{\max})) \big] .
$$
- second case: $M_{\max} - M_{\min}$ is even:
$$ h_{m}((n_{i})_{d},(\xi_{i})_{d}) = 
(-1)^{M_{\min}} \big[ (u(M_{\min}) - u(M_{\min+1})) + \cdots + (u(M_{\max}-1) - u(M_{\max})) + u(M_{\max}+1) \big] .
$$
By lemma \ref{second lemma alternative} and the hypothesis, $u$ is strictly decreasing. Moreover, $u$ is positive. In both cases above, $h_{m}((n_{i})_{d},(\xi_{i})_{d})$ is a sum of strictly positive terms, thus is strictly positive, thus is non-zero.
\end{proof}

Here is a concrete example in which the proposition \ref{prop alternative} applies : 

\begin{Corollary} \label{corollary alternative}Assume that $\xi_{2}=\ldots=\xi_{d}=1$ and $\xi_{1} = -1$. Then $h_{m}((n_{i})_{d};(\xi_{i})_{d}) \not= 0$.
\end{Corollary}

\begin{proof} By proposition \ref{prop alternative}, we only have to prove that $u$ is strictly decreasing on $[M_{\min},M_{\max}]$.
\newline We have $\displaystyle u(M) = \sum\limits_{M=m_{1}<\ldots<m_{d}<m} \frac{1}{m_{1}^{n_{1}} \ldots m_{d}^{n_{d}}}$, thus it is clear that $u$ is strictly decreasing.
\end{proof}

\section{Conditions on $l$ : from $\Lambda$-adjoint to adjoint $p$-adic cyclotomic multiple zeta values}

Until now, the main questions have been to identify the set of indices $(n_{i})_{d};(\xi_{i})_{d}$ such that
\newline  $\zeta_{p,\alpha}^{\Lambda\Ad}\big((n_{i})_{d};(\xi_{i})_{d}\big) \not= 0$, and the set of indices $(n_{i})_{d};(\xi_{i})_{d},m_{0},m$ such that $h_{m_{0},m}\big((n_{i})_{d};(\xi_{i})_{d}\big) \not= 0$.

Now we ask the following more precise question:

\begin{Question} Given an index $\big((n_{i})_{d};(\xi_{i})_{d}\big)$, 
what is the set of $l \in \mathbb{N}$ such that  $\zeta_{p,\alpha}^{\Ad}\big((n_{i})_{d};(\xi_{i})_{d};l\big)\not=0$? 
\end{Question}

This amounts to study the relation between the non-vanishing of $\Lambda$-adjoint $p$CMZV's and the non-vanishing of adjoint $p$CMZV's.

\subsection{Non-vanishing of adjoint vs $\Lambda$-adjoint $p$CMZVs, assuming the period conjecture}

The first type of information concerning the set $\{l \in \mathbb{N}_{\geq 0} \text{ }|\text{ }\zeta_{p,\alpha}^{\Ad}\big((n_{i})_{d};(\xi_{i})_{d};l\big)\not=0\}$ is about its size.

The period conjecture for $p$CMZV's (\cite{Deligne Goncharov} \S5.28, \cite{Yamashita}, conjectures 2.3.4, \cite{I-1}, conjecture 1.2.3) says that all polynomial equations in the $k_{N}$-algebra of $p$CMZV's (where $k_{N}$ is the $N$-th cyclotomic field) ``come from geometry''.

\begin{Proposition} \label{trichotomy}Let an index $\big((n_{i})_{d};(\xi_{i})_{d})$. Assume the period conjecture for $p$CMZV's. Then 

(i) The set 
$\{l \in \mathbb{N}_{\geq 0} \text{ }|\text{ }\zeta_{p,\alpha}^{\Ad}\big((n_{i})_{d};(\xi_{i})_{d};l\big)\not=0\}$ is either empty, or a singleton or infinite.

(ii) If we have $h_{p^{\alpha}}\big((n_{i})_{d};(\xi_{i})_{d}\big) \not=0$, then the set $\{l \in \mathbb{N}_{\geq 0} \text{ }|\text{ }\zeta_{p,\alpha}^{\Ad}\big((n_{i})_{d};(\xi_{i})_{d};l\big)\not=0\}$ is infinite.
\end{Proposition}

\begin{proof} (i) Let the alphabet $e_{0\cup\mu_{N}} = \{e_{x}\text{ }|\text{ }x= 0\text{ or }x^{N}=1\}$. To an index $\big((n_{i})_{d};(\xi_{i})_{d})$ of $p$CMZV's we associate the word $e_{0}^{n_{d}-1}e_{\xi_{d}} \cdots e_{0}^{n_{1}-1}e_{\xi_{1}}$. By their definition as iterated integrals, $p$-adic analogues of (\ref{eq:iterated integral CMZV}), $p$CMZV's satisfy the shuffle relation: for all words $w$, $w'$ on the alphabet $\{e_{x}\text{ }|\text{ }x= 0\text{ or }x^{N}=1\}$ we have $\zeta_{p,\alpha}(w)\zeta_{p,\alpha}(w') = \zeta_{p,\alpha}(w \sh w')$ where $\sh$ is the shuffle product on words on $e_{0\cup \mu_{N}}$. By the shuffle relation, applied to equation (\ref{eq:adjoint multiple zeta values}), each $\zeta_{p,\alpha}^{\Ad}\big((n_{i})_{d};(\xi_{i})_{d};l\big)$ is a $k_{N}$-linear combination of $p$CMZV's of weight $n_{1}+\cdots+n_{d}+l$.
	
The period conjecture for $p$CMZV's (\cite{Deligne Goncharov} \S5.28, \cite{I-1}, \S1) implies that the linear equations among $p$CMZV's are homogeneous for the weight. In other terms, the $k_{N}$-vector spaces generated by $p$CMZV's of a given weight are in direct sum. (Here, $k_{N}$ is regarded as the $k_{N}$-vector space generated by $p$CMZV's of weight $0$.)

Thus, if the set 
$\{l \in \mathbb{N}_{\geq 0} \text{ }|\text{ }\zeta_{p,\alpha}^{\Ad}\big((n_{i})_{d};(\xi_{i})_{d};l\big)\not=0\}$ is finite with at least two elements, equation (\ref{eq:theequation adjoint}) is a linear relation among $p$CMZV's which is not homogeneous for the weight, whence a contradiction.

(ii) The emptiness of this set would contradict the non-vanishing of $h_{p^{\alpha}}\big((n_{i})_{d};(\xi_{i})_{d}$ by equation  (\ref{eq:theequation adjoint}).
	
The equality between this set and a singleton would contradict the weight homogeneity of $p$CMZV's.
	
By (i), the only remaining possibility is that this set is infinite. 
\end{proof}

In particular, assuming the period conjecture for $p$CMZV's, the case where $\{l \in \mathbb{N}_{\geq 0} \text{ }|\text{ }\zeta_{p,\alpha}^{\Ad}\big((n_{i})_{d};(\xi_{i})_{d};l\big)\not=0\}$ is a singleton can happen only if $h_{p^{\alpha}}((n_{i})_{d};(\xi_{i})_{d}) = 0$. It corresponds to the case where equation (\ref{eq:theequation adjoint}) amounts to an equation of the form $\zeta_{p,\alpha}^{\Ad}((n_{i})_{d};(\xi_{i})_{d};l)=0$, i.e. the vanishing of a particular adjoint $p$CMZV. We think that this is rare. However, it can happen.

\subsection{The case of depth 1 and $N=1$}

As an illustration of the previous proposition, it can happen that $\{l \in \mathbb{N}_{\geq 0} \text{ }|\text{ }\zeta_{p,\alpha}^{\Ad}\big((n_{i})_{d};(\xi_{i})_{d};l\big)\not=0\}$ is infinite and at the same time its complement $\{l \in \mathbb{N}_{\geq 0} \text{ }|\text{ }\zeta_{p,\alpha}^{\Ad}\big((n_{i})_{d};(\xi_{i})_{d};l\big)=0\}$ is infinite too.

\begin{Example} Assume $N=1$ and $d=1$. We have, by the shuffle relation for $\Phi_{p,\alpha}$,
$$ \zeta_{p,\alpha}^{\Lambda\Ad}(n) = \sum_{l=1}^{\infty}  (-1)^{n} \zeta_{p,\alpha}(n+l)\Lambda^{l} = \sum_{l=1}^{\infty} \zeta^{\Ad}_{p,\alpha}(n+l)\Lambda^{l} . $$

We have $\zeta_{p,\alpha}(2n)=0$ for all $n\geq 1$, thus 
$\{l \in \mathbb{N}_{\geq 0}\text{ } |\text{ }\zeta_{p,\alpha}^{\Lambda\Ad} (n;l)=0\}$ is infinite.
We know that there are infinitely many $l\in \mathbb{N}$ such that $L_{p}(n+l,\omega^{-n-l})\not=0$ (see \S0.1), and it amounts to say, by equation (\ref{eq:zeta and L}), that the set $\{l \in \mathbb{N}_{\geq 0}\text{ } |\text{ }\zeta_{p,\alpha}^{\Lambda\Ad} (n;l)\not=0\}$ is infinite. 
Conjecturally, $\zeta_{p,\alpha}(2n+1)\not=0$ for all $n\geq 1$. 
\end{Example}

We note that we have used only a small part of the period conjecture, namely the homogeneity of the relations with respect to the weight.

\subsection{A bound on $l$}

Here is now a different type of information on the set $\{l \in \mathbb{N}_{\geq 0} \text{ }|\text{ }\zeta_{p,\alpha}^{\Ad}\big((n_{i})_{d};(\xi_{i})_{d};l\big)\not=0\}$. 
Assuming this set if non-empty, can we have a bound on its elements ?

More precisely : can we find an explicit upper bound on $\min \{l \in \mathbb{N}_{\geq 0} \text{ }|\text{ }\zeta_{p,\alpha}^{\Ad}\big((n_{i})_{d};(\xi_{i})_{d};l\big)\not=0\}$ ?

The answer is as follows. 

In \cite{I-1}, corollary 4.3.2, we prove that there exists constants $\kappa_{d},\kappa'_{d},\kappa''_{d}>0$, depending only on $d$, such that we have 
$v_{p}(\zeta_{p,\alpha}((n_{i})_{d};(\xi_{i})_{d})) \geq n - \kappa_{d} - \kappa'_{d}\log(n+\kappa'_{d})$ where $n=\sum_{i=1}^{d}n_{i}$. (It is possible using \cite{I-1} to obtain formulas for $\kappa_{d},\kappa'_{d},\kappa''_{d}$.) A more precise bound is obtained in \cite{Cha}. It is made explicit in the $N=1$ case : $\zeta_{p,-\infty}((n_{i})_{d}) \in \sum_{n\geq n_{0}} \frac{p^{n}}{n!}\mathbb{Z}_{p}$ (this has also been proved independently by Akagi, Hirose and Yasuda but unpublished).

Denoting by $v = v_{p} ((p^{\alpha})^{\sum_{i=1}^{d}n_{i}} \tilde{h}_{p^{\alpha}}((n_{i})_{d};(\xi_{i})_{d+1})$ (assuming this multiple harmonic sum is non-zero), equation (\ref{eq:theequation}) combined to lower bounds on the valuations of $p$CMZV's such as the ones above provides a formula for an upper bound on 
$\min \{l \in \mathbb{N}_{\geq 0} \text{ }|\text{ }\zeta_{p,\alpha}^{\Ad}\big((n_{i})_{d};(\xi_{i})_{d};l\big)\not=0\}$.

\subsection{Going beyond the hypothesis $p^{\alpha}>d$}

Until now, there has been an essential limitation in our method. The cyclotomic multiple harmonic sum $ $ can be non-zero only if $p^{\alpha}>d$, thus we have results of non-vanishing for adjoint $p$CMZV's $\zeta_{p,\alpha}^{\Ad}((n_{i})_{d};(\xi_{i})_{d};l)$ only when $p^{\alpha}>d$.

Here is how to deduce results when $p^{\alpha}\leq d$.

\begin{Proposition} Let $\alpha$ and $\alpha'$ be any elements of $\mathbb{Z} \cup \{\pm \infty\} - \{0\}$. 

(i) Assume that we have $\zeta_{p,\alpha}((n_{i})_{d};(\xi_{i})_{d})\not=0$. Then, there exists an index $((n'_{i})_{d'};(\xi'_{i})_{d'})$, with $\sum_{i=1}^{d'} n'_{i} = \sum_{i=1}^{d'}n_{i}$ and $d'\leq d$, such that $\zeta_{p,\alpha'}((n'_{i})_{d'};(\xi'_{i})_{d'})\not=0$.  

(ii) Assume that we have $\zeta_{p,\alpha}((n_{i})_{d};(\xi_{i})_{d})\not=0$. Then, there exists an index $((n'_{i})_{d'};(\xi'_{i})_{d'};l')$, with $\sum_{i=1}^{d'} n'_{i} + l' = \sum_{i=1}^{d'}n_{i}$ and $d'\leq d$, such that $\zeta^{\Ad}_{p,\alpha}((n'_{i})_{d'};(\xi'_{i})_{d'};l')\not=0$.  
(iii) Assume that we have $\zeta_{p,\alpha}^{\Ad}((n_{i})_{d};(\xi_{i})_{d};l)\not=0$. Then, there exists an index $((n'_{i})_{d'};(\xi'_{i})_{d'};l')$, with $\sum_{i=1}^{d'} n'_{i} + l' = \sum_{i=1}^{d'}n_{i}+l$ and $d'\leq d$, such that $\zeta_{p,\alpha}^{\Ad}((n'_{i})_{d'};(\xi'_{i})_{d'};l')\not=0$.
\end{Proposition}

\begin{proof} (i) This follows from \cite{I-3}, corollary 2.2.2.

(ii) This follows from \cite{II-1}, corollary 2.2.5.
	
(iii) This follows from \cite{I-3}, corollary 2.2.2 and \cite{II-1}, corollary 1.2.5 combined.
\end{proof}

Combining proposition 5.7 with results of \S1 to \S4, we deduce non-vanishing results for  $\zeta^{\Ad}_{p,\alpha}((n_{i})_{d};(\xi_{i})_{d};l)$ or  $\zeta_{p,\alpha}((n_{i})_{d};(\xi_{i})_{d})$ with $p^{\alpha} \leq d$. This would not have been possible if we had not generalized the definition of $p$CMZV's in \cite{I-1} by introducing the parameter $\alpha$. We would have had to restrict the non-vanishing results to $p>d$, or (depending on the convention) $q>d$ where $q$ is such that $\mathbb{F}_{q}$ contains a primitive $N$-th root of unity. We also note that the case $\alpha = - \infty$ corresponds to Coleman integration. 

In parallel to this proof, for any algebraic relation among (adjoint) $p$CMZV's ($\zeta_{p,\alpha}$ for a fixed $\alpha$, or $\zeta^{\Ad}_{p,\alpha}$ for a fixed $\alpha$), if one of the terms does not vanish, then there exists at least another non-zero term, since the sum of all terms is zero. This can be used to obtain the existence of certain words $w$ such that $\zeta_{p,\alpha}(w)\not=0$. However, these words $w$ are not explicit, unlike those found in the sections \S1 to \S4.  We have used it in \S2.5. Here is another example :

\begin{Example} $p$CMZV's satisfy the double shuffle relations : \cite{Besser Furusho, Furusho Jafari} for $N=1$ and by \cite{Yamashita} for any $N$ (this is for particular values of $\alpha$, and this implies that this is true for all $\alpha$ as explained in \cite{II-1}).
	
There is a long list of examples of equations obtained as consequences of the double shuffle relations. For example, the ``sum formula'' \cite{Granville, Hoffman}: for all $n \in \mathbb{N}^{\ast}$, $n \geq 2$, and for all $d \in \mathbb{N}^{\ast}$, we have $\zeta_{p,\alpha}(n) = 
\sum\limits_{\substack{n_{1},\ldots,n_{d} \in \mathbb{N}^{\ast} \\ n_{d} \geq 2}}\zeta_{p,\alpha}(n_{1},\ldots,n_{d})$.
As explained in the introduction, for any $(p,\alpha)$, there are infinitely many $n$'s such that $\zeta_{p,\alpha}(n) \not=0$. As a consequence, there exists a sequence of words $(w_{l})_{l\in \mathbb{N}}$  such that $\weight(w_{l}) \underset{l\rightarrow\infty}{\rightarrow} \infty$ and $\depth(w_{l})=d$, and $\zeta_{p,\alpha}(w_{l}) \not= 0$ for all $l \in \mathbb{N}$.
\end{Example}


\begin{thebibliography}{50}
\bibitem[B]{Balliet} K. Balliet - \emph{On the prime numbers in intervals}, master thesis, 2015, arXiv:1706.01009.
\bibitem[BF]{Besser Furusho} A. Besser, H. Furusho - \emph{The double shuffle relations for p-adic multiple zeta values}, AMS Contemporary Math., 416 (2006) pp. 9-29.
\bibitem[Cha]{Cha} A. Chatzitamatiou - \emph{On integrality of $p$-adic iterated integrals,} J. Algebra, 474, n°15 (March 2017), pp. 240-270.
\bibitem[Ch]{Chen} K. T. Chen - \emph{Iterated path integrals}, Bull. Amer. Math. Soc, Vol. 83, n°5 (1977) pp. 831-879.
\bibitem[Col]{Coleman} R. Coleman - \emph{Dilogarithms, regulators and $p$-adic $L$-functions}, Invent. Math., 69, n°2 (1982), pp. 171-208.
\bibitem[D]{Deligne} P. Deligne - \emph{Le groupe fondamental de la droite projective moins trois points}, Galois Groups over $\mathbb{Q}$ (Berkeley, CA, 1987), Math. Sci. Res. Inst. Publ. 16, Springer-Verlag, New York, 1989.
\bibitem[DG]{Deligne Goncharov} P. Deligne, A. B. Goncharov - \emph{Groupes fondamentaux motiviques de Tate mixtes}, Ann. Sci. Ecole Norm. Sup. 38, n°1 (2005), pp. 1-56
\bibitem[F1]{Furusho 1} H. Furusho - \emph{p-adic multiple zeta values I -- p-adic multiple polylogarithms and the p-adic KZ equation}, Invent. Math., 155, (2004) n°2, pp. 253-286.
\bibitem[F2]{Furusho 2} H. Furusho - \emph{p-adic multiple zeta values II -- tannakian interpretations}, Amer.J.Math, Vol 129, n°4, (2007), pp. 1105-1144.
\bibitem[FJ]{Furusho Jafari} H. Furusho, A. Jafari - \emph{Regularization and generalized double shuffle relations for p-adic multiple zeta values}, Compositio Math. 143 (2007), pp. 1089-1107.
\bibitem[Gr]{Granville} A. Granville - \emph{A decomposition of Riemann's zeta-function}, in Analytic number theory, London Math. Soc. Lecture Note Ser., vol. 247, Cambridge Univ. Press, (1997), pp. 95-102. 
\bibitem[J I-1]{I-1} D. Jarossay - \emph{A bound on the norm of overconvergent $p$-adic multiple polylogarithms}, 
J. Number Theory, Vol. 205 (December 2019), pp. 81-121.
\bibitem[J I-2]{I-2} D. Jarossay - \emph{Pro-unipotent harmonic actions and a computation of $p$-adic cyclotomic multiple zeta values}, arXiv:1501.04893, submitted.
\bibitem[J I-3]{I-3} D. Jarossay - \emph{Pro-unipotent harmonic actions and dynamical properties of $p$-adic cyclotomic multiple zeta values}, Alg. Number Th. 14 (2020) pp. 1711-1746.
\bibitem[J II-1]{II-1} D. Jarossay - \emph{Adjoint cyclotomic multiple zeta values and cyclotomic multiple harmonic values}, arXiv:1412.5099, submitted.
\bibitem[J II-3]{II-3} D. Jarossay - \emph{Cyclotomic multiple harmonic values regarded as periods}, arxiv:1601.01159.
\bibitem[J III-1]{III-1} D. Jarossay - \emph{$p$-adic cyclotomic multiple zeta values for roots of unity of order divisible by $p$}, arxiv:1708.08009.
\bibitem[H]{Hoffman} M. E. Hoffman - \emph{Multiple harmonic series}, Pacific J. Math. 152 (1992), no. 2, pp. 275-290.
\bibitem[KL]{KL} H. Kadiri, A. Lumley - \emph{Short effective intervals containing primes}, arXiv:1407.7902.
\bibitem[RS]{RS} O. Ramaré, Y. Saouter - \emph{Short effective intervals containing primes}, J. Number Theory 98 (2003), pp. 10-33.
\bibitem[So]{Soule} C. Soul\'{e} - \emph{On higher p-adic regulators} Algebraic K-theory, Evanston 1980 (Proc. Conf., Northwestern Univ., Evanston, Ill., 1980), pp. 372-401. Lect. Notes Math. 854. Berlin, New York, Springer 1981.
\bibitem[U1]{Unver 1} S. Unver - \emph{$p$-adic multi-zeta values} J. Number Theory n°108 (2004), pp. 111-156 .
\bibitem[U2]{Unver 2} S. Unver - \emph{Cyclotomic $p$-adic multi-zeta values in depth two}, Manuscripta Math. 149, (2016) no. 3-4, pp. 405-441.
\bibitem[W]{Washington} L. C. Washington - \emph{$p$-adic $L$-functions and sums of powers}, J. Number Theory 69 (1998), pp. 50-61.
\bibitem[Yam]{Yamashita} G. Yamashita - \emph{Bounds for the dimension of $p$-adic multiple $L$-values spaces}, Documenta Mathematica, Extra Volume Suslin (2010) pp. 687-723.
\end{thebibliography}
\end{document}